\documentclass[12pt]{article}

\usepackage{amsmath}
\usepackage{latexsym}
\usepackage{amssymb}
\usepackage[mathcal]{euscript}

\def\ee{\epsilon}
\def\eps{\epsilon}
\def\aa{f} 
\def\vrln#1{#1}

\newcommand\bbR{{\mathbb R}}
\newcommand{\bbN}{{\mathbb N}}

\newtheorem{theorem}{Theorem}

\newtheorem{lemma}[theorem]{Lemma}

\usepackage{graphicx}
\usepackage{color}

\begin{document}

\title{Curvature line parametrization from circle patterns}

\author{A.I.~Bobenko
 and S.P.~Tsarev\thanks{On leave from: Krasnoyarsk State Pedagogical
University, Russia. SPT acknowledges partial financial support
from the DFG Research Unit 565 `Polyhedral Surfaces'' (TU-Berlin) and
the RFBR grant 
 06-01-00814.}}

   \maketitle
\begin{center}
Department of Mathematics \\
 Technische Universit\"at Berlin \\
Berlin, Germany \\[1ex] 
e-mails: \\
\texttt{bobenko@math.tu-berlin.de}\\
\texttt{tsarev@math.tu-berlin.de} \ 
\texttt{tsarev@newmail.ru}\\
\end{center}

\medskip

\begin{abstract}
 We study local and global approximations of smooth nets of curvature lines and
smooth conjugate nets  by respective discrete nets (circular nets and planar quadrilateral nets)
with edges of order $\epsilon$.  It is shown that choosing the
points of discrete nets on the smooth surface one can obtain the
order $\epsilon^2$ approximation globally. Also a simple geometric construction for
approximate determination of principal directions of smooth surfaces is given.
\end{abstract}

\section{Introduction}

Discrete conjugate nets, defined as mappings ${\bf Z}^2
\longrightarrow \bbR^3$ with the condition that each elementary
quadrangle is flat, and discrete nets of curvature lines, defined as
discrete conjugate nets with additional property of circularity of
the four vertices of every elementary quadrangle, play an important
role in the contemporary discrete differential geometry (see
e.g.~\cite{BS}
for a review) and have important 
applications to computer-aided geometric design \cite{BPot}.

This paper is devoted to a carefull study of the order of approximation 
of discrete
nets to a given smooth conjugate net or the net of curvature lines
on a smooth surface that can be achieved. We also present a remarkably simple 
approximate construction of principal directions on a given smooth surface using
the standard primitive of discrete differential geometry: an elementary circular
quadrangle inscribed into the surface. Roughly speaking, our
results show that the previously known upper bounds  (\cite{BMS})
can be made one order better for conjugate nets and discrete nets of
curvature lines. Moreover, one can impose one additional geometric
requirement: all vertices of approximating discrete nets should lie
on the original smooth surface.

The paper is organized as follows. In Section~\ref{sec2} we present a method
of construction of principal directions in a neighborhood of a non-umbilic point
on a smooth surface using only the four points of intersection of the surface with an
infinitesimal circle. Our method is in fact based on a third-order approximation
of the given smooth surface and gives an $\ee^2$-approximation of principal
directions if the infinitesimal circle has the radius $\ee$. Different methods
for approximation of principal directions have been proposed recently, see for example
\cite{GI}. Re-meshing techniques based on determination of principal directions
are widely used in contemporary computer graphics (\cite{ACDLM}).

Section~\ref{sec-local} is devoted to the study of the order of approximation of 
a single elementary quad of a discrete conjugate net (resp.\ discrete
circular net) to the respective infinitesimal quadrangle of a the smooth net of conjugate lines (resp.\ curvature lines)
on the original smooth surface. The results of this Section are used in 
Section~\ref{sec-global} where we prove the main global results on the order
of approximation of respective discrete nets to a given smooth conjugate net 
or the net of curvature lines on a smooth surface.

We use the following terminology and notations: $d(A,S)$ is the
distance from a point $A$ to another point $S$ (or a line, plane,
etc.). A point $A$ is said to be $\epsilon^k$-close to $S$ if
$d(A,S) =O(\epsilon^k)$, that is $d(A,S) \leq C \epsilon^k$  for
$\epsilon \to 0$.
$A$ is said to lie at $\epsilon^k$-distance from $S$, if $d(A,S)
\sim C \epsilon^k$ for some $C \neq 0$, we will denote this
hereafter as $d(A,S) \simeq \epsilon^k$. The same terminology and
notations will be used in characterization of the asymptotic
behavior for other (numeric) quantities depending on $\ee$.

\section{Principal curvature directions from circles}\label{sec2}

In this Section we show that one can use the elementary quadrangle
of a discrete circular net inscribed into a smooth surface (see
Section~\ref{sec32} for more details) for a rather precise
construction of the principal directions on a smooth surface.

\begin{itemize}
 \item  \textbf{\underline{Euclidean construction.}} Take three points  $A$, $B$,  $C$ on a smooth surface
$\aa:\Omega \longrightarrow {\bbR}^3$ without umbilic points, such
that $d(A,B)\simeq \ee$, $d(B,C)\simeq \ee$ and the angle $ \angle
ABC$ is not $\ee$-close to $0$ or $2\pi$.Then the circle $\omega$
passing through   $A$, $B$,  $C$ has radius of order $O(\ee)$ as
well; from elementary topological considerations one concludes that
$\omega$ has at least one more point of intersection with
$\aa(\Omega)$. Denote it (or one of them) as $D$. Take the
straight lines $(AC)$,  $(BD)$ and the point of their intersection
$Z=(AC)\cap(BD)$ (see Fig.~\ref{fig1}). The bisectors $b_1$ and
$b_2$ of the angles formed at $Z$ by $(AC)$ and  $(BD)$ are
obviously perpendicular. In fact $b_1$ and $b_2$ approximate the
directions of smooth curvature lines on the surface: as we prove
below, the directions of  $b_1$ and $b_2$ are $\ee$-close to the
principal directions at any point $M\in\aa(\Omega)$ which is
$\ee$-close to $A$, $B$,  $C$, $D$.
 In fact one can find the principal directions with only $O(\ee^2)$-error at
some specific point $Z_\aa$ of $\aa(\Omega)$ which ``lies over
$Z$''. Note that the order of the points $A$, $B$,  $C$, $D$ on the circle $\omega$
is important for finding the principal direction with the error of order $\ee^2$.

Namely, let $Z_\aa$ be the intersection of  $\aa(\Omega)$ with
the straight line $l$ passing through $Z$ and perpendicular to the
plane  $\pi_\epsilon=(ABC)$. Take the orthogonal (parallel to $l$)
projections  $p_1$ and $p_2$ of the bisectors $b_1$ and $b_2$ onto
the tangent plane $\pi_{\aa}$  of $\aa(\Omega)$  at  $Z_\aa$.
\end{itemize}

\begin{theorem} \label{Th-pd}
For a general smooth surface $\aa(\Omega)$ in a neighborhood of a non-umbilic point 
the Euclidean construction given above produces the directions
$p_1$, $p_2$ which are $\ee^2$-close to the exact principal
directions at $Z_\aa$.
\end{theorem}

In order to prove this Theorem we establish at first some simple but
remarkable facts about planar quadrics and Dupin cyclides
(Lemma~\ref{thd} and Theorem~\ref{th-dupen}) and give another
(M\"obius-invariant) construction of the principal directions
(Fig.~\ref{fig2}). In fact, the intersection curve
$\delta=\aa(\Omega) \cap \pi_\epsilon$ is $\epsilon^2$-close to a
quadric --- the Dupin indicatrix of $\aa(\Omega)$ at $P$ 
(see the Appendix); 
this explains in particular why we may
assume that the number of intersection points of  $\omega$ and
$\aa(\Omega)$ is \textit{exactly} four. 
If we take then an {\em
arbitrary} circle $\omega$ in a plane, intersecting  a
quadric in points $A$, $B$, $C$, $D$ (cf.\ Fig.~\ref{fig1} for the elliptic case; the order of
the points  $A$, $B$, $C$, $D$ on the circle is not necessarily consecutive for Lemma~\ref{thd} below
to hold), 
then the directions of the angles formed by the straight lines 
$(AC)$ and $(BD)$ are parallel to axes of the quadric:

\begin{figure}[htbp]
\begin{center}
\begin{picture}(0,0)%
\includegraphics{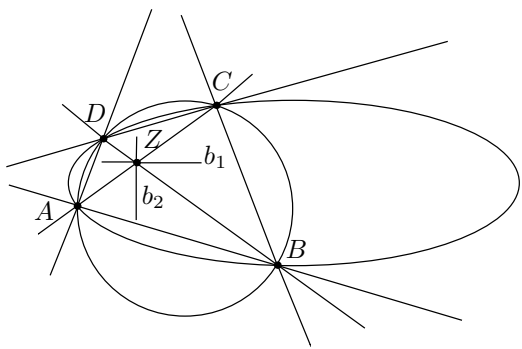}%
\end{picture}%
\setlength{\unitlength}{4144sp}%
\begingroup\makeatletter\ifx\SetFigFont\undefined%
\gdef\SetFigFont#1#2#3#4#5{%
  \reset@font\fontsize{#1}{#2pt}%
  \fontfamily{#3}\fontseries{#4}\fontshape{#5}%
  \selectfont}%
\fi\endgroup%
\begin{picture}(3089,2032)(46,-1337)
\put(878,-146){\makebox(0,0)[lb]{\smash{{\SetFigFont{10}{12.0}{\familydefault}{\mddefault}{\updefault}{\color[rgb]{0,0,0}$Z$}%
}}}}
\put(870,-483){\makebox(0,0)[lb]{\smash{{\SetFigFont{10}{12.0}{\familydefault}{\mddefault}{\updefault}{\color[rgb]{0,0,0}$b_2$}%
}}}}
\put(1242,-241){\makebox(0,0)[lb]{\smash{{\SetFigFont{10}{12.0}{\familydefault}{\mddefault}{\updefault}{\color[rgb]{0,0,0}$b_1$}%
}}}}
\put(1730,-802){\makebox(0,0)[lb]{\smash{{\SetFigFont{10}{12.0}{\familydefault}{\mddefault}{\updefault}{\color[rgb]{0,0,0}$B$}%
}}}}
\put(1288,199){\makebox(0,0)[lb]{\smash{{\SetFigFont{10}{12.0}{\familydefault}{\mddefault}{\updefault}{\color[rgb]{0,0,0}$C$}%
}}}}
\put(522,  9){\makebox(0,0)[lb]{\smash{{\SetFigFont{10}{12.0}{\familydefault}{\mddefault}{\updefault}{\color[rgb]{0,0,0}$D$}%
}}}}
\put(225,-576){\makebox(0,0)[lb]{\smash{{\SetFigFont{10}{12.0}{\familydefault}{\mddefault}{\updefault}{\color[rgb]{0,0,0}$A$}%
}}}}
\end{picture}%
\end{center}
\caption{Construction of the axes directions of a quadric: $A$, $B$,
$C$, $D$ are the four points of intersection with a circle, $b_1$
and $b_2$ are the bisectors of the angles between $(AC)$ and
$(BD)$. The bisectors $b_1$ and $b_2$ are parallel to the axes of
the quadric} \label{fig1}
\end{figure}

\begin{lemma}\label{thd} For any 
plane quadric with four point of intersection with a circle 
the bisectors  $b_1$ and $b_2$
(shown on Fig.~\ref{fig1} for the elliptic case) are parallel to the
axes of the quadric.
\end{lemma}
\textsl{Proof.} For simplicity we give the proof only for the case
of a non-degenerate central quadric
$Q(x,y)=(K_1x^2+ K_2y^2) - 2\ee^2=0$; the other cases can be proved along the same lines.
The coordinates of the intersection points $A$, $B$,
$C$, $D$ satisfy the system
\begin{equation}
\left\{\begin{array}{l}\label{Qomega}
Q(x,y)=(K_1x^2+ K_2y^2) - 2\ee^2=0, \\
\omega(x,y) = (x-x_\ee)^2 + (y-y_\ee)^2  - r^2_\ee =0.
\end{array}\right.
\end{equation}
They also lie on the one-parametric sheaf of quadrics $Q_\lambda$:
$\{Q(x,y)+\lambda \omega(x,y)=0\}$. Obviously
the directions of the axes 
of all $Q_\lambda$ are the same. On the other hand, for
some values of $\lambda$ this quadric is degenerate. These values are the three roots $\lambda_i$ 
of the equation
\begin{equation}
\left\vert \begin{array}{ccc}
K_1+ \lambda & 0 & -2\lambda x_\ee \\
0 & K_2+\lambda &  -2\lambda y_\ee \\
 -2\lambda x_\ee  &  -2\lambda y_\ee  & 2\ee^2 - \lambda r^2_\ee
\end{array} \right\vert  = 0.
\end{equation}
Two of them correspond to the degenerate quadrics composed of the
pairs of straight lines $\{(AB),(CD)\}$ and $\{(AD),(BC)\}$ shown on
Fig.~\ref{fig1}. The third degenerate quadric is the pair of
straight lines 
$\{(AC),(BD)\}$, shown on Fig.~\ref{fig1}.

Since the bisectors of these pairs of lines define exactly the axes 
of the respective degenerate quadric $Q_\lambda$, we have proved
that their directions are the same as for $Q_0=\{Q(x,y)=0\}$. 
$\Box$

Note that we have not used the ordering of the points $A$, $B$,
$C$, $D$ on the circle $\omega$ or on the quadric, so the statement of Lemma~\ref{thd} holds
true also for the bisectors of the angles formed by the lines $(AB)$ and $(CD)$,
as well as $(AD)$ and $(BC)$. On the other hand in order to obtain the approximation
of order $\ee^2$ below we \textit{fix} the ordering as shown on Fig.~\ref{fig1}.

Next we need to reformulate the Euclidean construction given above
in terms of M\"obius geometry:

\begin{itemize}
 \item \textbf{\underline{M\"obius-invariant construction}. }
We start with a smooth surface $\aa: \Omega \longrightarrow \bbR^3$.
Take a sphere $S$ tangent to $\aa(\Omega)$ at some point $P$,
such that for its radius $R$ and the principal curvatures $K_1$,
$K_2$ ($K_1 <K_2$) of $\aa(\Omega)$ at $P$ the inequalities $K_1 < 1/R <
K_2$ hold (with non-infinitesimal differences $|K_1 - 1/R|$ and $|
1/R- K_2|$ ). In this case $S$ intersects an $\ee$-patch
$\Omega_\ee$ on $\aa(\Omega)$ along two smooth lines $\gamma_1$,
$\gamma_2$ (Fig.~\ref{fig2}): $\Omega_\aa \cap S = \gamma_1 \cup
\gamma_2$, the angle between  $\gamma_1$ and $\gamma_2$ is non-zero
and  non-infinitesimal (this follows from the classical lemma of
Morse, cf.\ for example \cite{Mi}).
Take any 
circle $\omega \subset S$ of radius $\simeq \ee$ with $P$ inside
such that the distances between $P$ and the four points of
intersection of $\omega$ and $\Omega_\ee$: $\omega \cap \gamma_1 =A
\cup C$,   $\omega \cap \gamma_2 =B \cup D$ are of order $O(\ee)$ as
well. Construct two auxiliary circles $\omega_1$, $\omega_2$ on $S$
defined by the triples of points $\{A,P,C\}$ and $\{B,P,D\}$
respectively (not shown on Fig.~\ref{fig2}, $\omega_i$ are very
close to $\gamma_i$). Then the bisectors  $b_1$, $b_2$  of the
angles formed by $\omega_1$, $\omega_2$  at $P$ obviously lie in the
tangent plane to $\aa(\Omega)$ at $P$. As we show below, $b_1$ and $b_2$
are $\ee^2$-close to the principal directions of $\aa(\Omega)$ at $P$. 
\end{itemize}

\begin{figure}[htbp]
\begin{center}
\begin{picture}(0,0)
\includegraphics{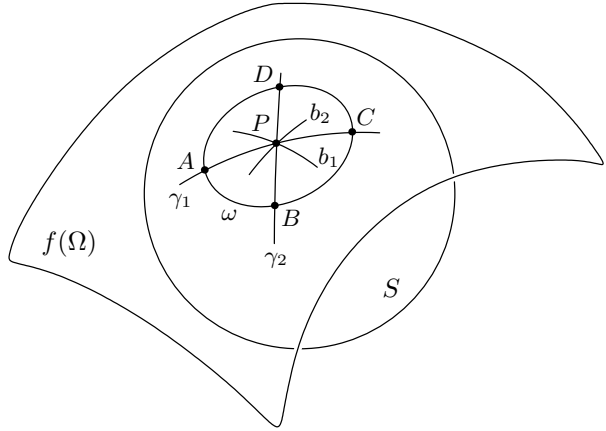}%
\end{picture}%
\setlength{\unitlength}{4144sp}%
\begingroup\makeatletter\ifx\SetFigFont\undefined%
\gdef\SetFigFont#1#2#3#4#5{%
  \reset@font\fontsize{#1}{#2pt}%
  \fontfamily{#3}\fontseries{#4}\fontshape{#5}%
  \selectfont}%
\fi\endgroup%
\begin{picture}(3578,2559)(383,-1906)
\put(1852,173){\makebox(0,0)[lb]{\smash{{\SetFigFont{9}{10.8}{\familydefault}{\mddefault}{\updefault}{\color[rgb]{0,0,0}$D$}%
}}}}
\put(2473,-89){\makebox(0,0)[lb]{\smash{{\SetFigFont{9}{10.8}{\familydefault}{\mddefault}{\updefault}{\color[rgb]{0,0,0}$C$}%
}}}}
\put(1842,-123){\makebox(0,0)[lb]{\smash{{\SetFigFont{9}{10.8}{\familydefault}{\mddefault}{\updefault}{\color[rgb]{0,0,0}$P$}%
}}}}
\put(1396,-352){\makebox(0,0)[lb]{\smash{{\SetFigFont{9}{10.8}{\familydefault}{\mddefault}{\updefault}{\color[rgb]{0,0,0}$A$}%
}}}}
\put(2629,-1118){\makebox(0,0)[lb]{\smash{{\SetFigFont{10}{12.0}{\familydefault}{\mddefault}{\updefault}{\color[rgb]{0,0,0}$S$}%
}}}}
\put(588,-836){\makebox(0,0)[lb]{\smash{{\SetFigFont{10}{12.0}{\familydefault}{\mddefault}{\updefault}{\color[rgb]{0,0,0}$f(\Omega)$}%
}}}}
\put(2021,-692){\makebox(0,0)[lb]{\smash{{\SetFigFont{9}{10.8}{\familydefault}{\mddefault}{\updefault}{\color[rgb]{0,0,0}$B$}%
}}}}
\put(2193,-58){\makebox(0,0)[lb]{\smash{{\SetFigFont{9}{10.8}{\familydefault}{\mddefault}{\updefault}{\color[rgb]{0,0,0}$b_2$}%
}}}}
\put(2244,-322){\makebox(0,0)[lb]{\smash{{\SetFigFont{9}{10.8}{\familydefault}{\mddefault}{\updefault}{\color[rgb]{0,0,0}$b_1$}%
}}}}
\put(1665,-657){\makebox(0,0)[lb]{\smash{{\SetFigFont{9}{10.8}{\familydefault}{\mddefault}{\updefault}{\color[rgb]{0,0,0}$\omega$}%
}}}}
\put(1351,-562){\makebox(0,0)[lb]{\smash{{\SetFigFont{9}{10.8}{\familydefault}{\mddefault}{\updefault}{\color[rgb]{0,0,0}$\gamma_1$}%
}}}}
\put(1921,-912){\makebox(0,0)[lb]{\smash{{\SetFigFont{9}{10.8}{\familydefault}{\mddefault}{\updefault}{\color[rgb]{0,0,0}$\gamma_2$}%
}}}}
\end{picture}%
\end{center}
\caption{Approximate M\"obius-invariant construction of the
principal directions. A sphere $S$ is tangent to $\aa(\Omega)$ at
$P$, it intersects $\aa(\Omega)$ along two smooth lines
$\gamma_1$, $\gamma_2$. A circle $\omega \subset S$ of radius
$\simeq \ee$ with $P$ inside gives four points: $\omega \cap
\gamma_1 =A \cup C$,   $\omega \cap \gamma_2 =B \cup D$,
the points $A$, $B$, $C$, $D$ should lie on $\ee$-distance from $P$.
The bisectors  $b_1$, $b_2$  of the angles formed by the circles
$\omega_1$, $\omega_2$ (not shown) passing through the triples of
points $\{A,P,C\}$ and $\{B,P,D\}$  are $\ee^2$-close to  the principal
directions of $\aa(\Omega)$ at $P$}\label{fig2}
\end{figure}

First we prove that for Dupin cyclides (standard M\"obius
primitives, playing the role of the osculating paraboloids of the
classical Euclidean differential geometry, see \cite{Da1,HJ})  $b_i$
are \emph{exactly} the principal directions:
\begin{theorem} \label{th-dupen}
For a smooth non-umbilic 
point $P$ on a Dupin cyclide $\mathcal{D}$, the bisectors  $b_1$ and
$b_2$
of the M\"obius-invariant construction give the principal directions
of  $\mathcal{D}$ at $P$.
\end{theorem}
\textsl{Proof.} Since our construction of  $\omega_i$ is
M\"obius-invariant, we may transform $\mathcal{D}$ into one of the
normalized Dupin cyclides: a torus, a circular cone or
a circular cylinder. We give below a detailed proof for the case of
a torus; the other cases can be easily proved along the same lines.

Making, if necessary, another M\"obius transformation, we can reduce
the configuration to the following: the torus is given by the
equation
\begin{equation}\label{tor1}
(R^2 + x^2 + y^2 + z^2 - r^2)^2 = 4R^2(x^2 + y^2),
\end{equation}
the point $P$ has the coordinates $(R+r,0,0)$ and the sphere $S$ is
given by the equation
\begin{equation}\label{sf1}
\big(x+ \rho -(R+r)\big)^2 + y^2 + z^2 = \rho^2.
\end{equation}
After a M\"obius inversion $(x,y,z) \rightarrow (x_1,y_1,z_1)$ with
the center $P$:
\begin{equation}\label{mob}
\left\{
\begin{array}{l}
x_1= \displaystyle \frac{x-R-r}{ (x- (R+r) )^2 + y^2 + z^2},\\[1.3em]
y_1= \displaystyle \frac{y}{ (x- (R+r) )^2 + y^2 + z^2},\\[1.3em]
z_1= \displaystyle \frac{z}{ (x- (R+r) )^2 + y^2 + z^2},
\end{array}
\right.
\end{equation}
the sphere (\ref{sf1}) will be transformed into a plane
$x_1=\mathrm{const}=\mu$. We prove now that the intersection of the
image of the torus with any plane $x_1=\mu$ is a quadric of the form
$ay_1^2+ bz_1^2 =\kappa$. In fact, from (\ref{mob}) we have:
\begin{equation}\label{x1mu}
\Big(\big(x- (R+r)\big)^2 + y^2 + z^2\Big) = (x-R-r)/\mu
\end{equation}
so
\begin{equation}\label{z2}
z^2 = (x-R-r)/\mu - \Big(\big(x- (R+r)\big)^2 + y^2 \Big).
\end{equation}
(\ref{tor1}) may be now reshaped to
$$
\Big(\big(x- (R+r)\big)/\mu +2(R+r)x -2rR\Big)^2 = 4R^2(x^2 + y^2),
$$
so one can express $y^2$ in terms of $x$, $R$, $r$, $\mu$.
Substituting (\ref{x1mu}), (\ref{z2}) and the found expression for
$y^2$ into $ay_1^2 + bz_1^2 = (ay^2 + bz^2)\Big/\Big(\big(x-
(R+r)\big)^2 + y^2 + z^2\Big)$ one gets a rational  expression
containing only $x$, the constants  $R$, $r$, $\mu$ and $a$, $b$. A
straightforward calculation shows that for $a= r(1+2\mu(R+r))$,
$b=(1+2r\mu)(r+R)$ this expression is a constant $\kappa=-R(1+2\mu
R)(1+2\mu r)(1+2\mu r + 2\mu R)/4$.

Since the images of the circles $\omega_1$, $\omega_2$ defined above
are now straight lines passing through the points of intersection of
the image of the circle $\omega$ and the found quadric on the plane
$x_1=\mu$, the statement of this Theorem is now equivalent to the
statement of Lemma~\ref{thd}. $\Box$

\begin{theorem} \label{Th-Mic}
For an arbitrary smooth surface $\aa(\Omega)$ in a neighborhood of a
non-umbilic point $P$ the M\"obius-invariant construction produces
the directions $b_1$, $b_2$ which are $\ee^2$-close to the exact
principal directions at $P$.
\end{theorem}
\textsl{Proof.} First we perform a M\"obius transformation that maps
the sphere $S$ and the circle $\omega$ into themselves and brings
$P$ into the (Euclidean) center of $\omega$ on $S$, i.e.\ into the
point $P' \in S$ such that the distance from $P'$ to all the points
of $\omega$ are equal. One can easily check that our requirement
$d(P,A) \simeq\ee$, \ldots ,  $d(P,D) \simeq\ee$ guarantees that
after the transformation these distances will remain of order
$\simeq\ee$. Inside the circle  $\omega$ on $S$ and in its
$\ee$-neighborhood in $\bbR^3$ where $\aa(\Omega_\ee)$ lies, the
Jacobian of this transformation is limited so the distances of order
$O(\ee^p)$, $p \geq 2$ will be transformed into distances of the
same order. Since the angles do not change, $\ee^2$-close directions
at $P$ remain $\ee^2$-close at $P'$  and vice versa. Introduce now a
Cartesian coordinate system such that $P'$ is its origin, the
coordinate plane $(xy)$ is tangent to $S$ (so to $\aa(\Omega)$ as well) 
and the $x,y$-axes are the principal directions of $\aa(\Omega)$ at $P'$. 
The plane where $\omega$ lies will be given by the equation
$z=c_1\ee^2$, the circle $\omega$ will be defined as $x^2 +
y^2=c^2_2\ee^2$ and $\aa(\Omega)$ will be defined in a neighborhood of
$P'$ by its Taylor expansion
\begin{equation}\label{eq-urp}
z= (K_1 x^2 + K_2 y^2)/2 + (a_{30}x^3 + a_{21}x^2 y + a_{12}x y^2 +
a_{03}y^3) + o(\ee^3).
\end{equation}
So the intersection points $P_1=A$, $P_2=B$, $P_3=C$, $P_4=D$ of
$\omega$ and $\aa(\Omega)$  will have respectively the coordinates 
$x_i=\ee x_{i0} + \ee^2 x_{i1} + o(\ee^2)$ with $(x_{i0},y_{i0})$
being the solutions of
\begin{equation}\label{xy1}
\left\{ 
\begin{array}{l}
x^2 + y^2 = c^2_2,\\
(K_1 x^2 + K_2 y^2) = 2c_1,
\end{array}
\right.
\end{equation}
so we see that  
\begin{equation}\label{plusminus}
(x_{10},y_{10})=(x_{20},-y_{20})=(-x_{30},-y_{30})=(-x_{40},y_{40}).
\end{equation}
The next terms $(x_{i1},y_{i1})$ are found after substitution of
$(x_{i},y_{i})$ into
\begin{equation}\label{xy}
\left\{
\begin{array}{l}
x^2 + y^2 = c^2_2,\\
(K_1 x^2 + K_2 y^2)/2 + \ee (a_{30}x^3 + a_{21}x^2 y + a_{12}x y^2 +
a_{03}y^3) = c_1,
\end{array}
\right.
\end{equation}
cancellation of the second-order terms using (\ref{xy1}) and
retaining only third-order terms:
\begin{equation}\label{xy2}
\left\{
\begin{array}{l}
2x_{i0}x_{i1} + 2y_{i0}y_{i1} = 0,\\
K_1 x_{i0}x_{i1} + K_2 y_{i0}y_{i1} + (a_{30}x_{i0}^3 +
a_{21}x_{i0}^2 y_{i0} + a_{12}x_{i0} y_{i0}^2 + a_{03}y_{i0}^3) = 0.
\end{array}
\right.
\end{equation}

From (\ref{xy2}) and (\ref{plusminus}) we easily conclude that
\begin{equation}\label{f12}
(x_{11},y_{11})=(x_{31},y_{31}), (x_{21},y_{21})=(x_{41},y_{41}),
\end{equation}
so the $O(\ee^2)$-shifted lines $l_{11}$, $l_{21}$ defined by the
pairs of points $\{(\ee x_{10} + \ee^2 x_{11}, \ee y_{10} + \ee^2
y_{11} ),(\ee x_{30} + \ee^2 x_{31},
 \ee y_{30} + \ee^2 y_{31} )\}$ and
$\{(\ee x_{20} + \ee^2 x_{21}, \ee y_{20} + \ee^2 y_{21} ),(\ee
x_{40} + \ee^2 x_{41},  \ee y_{40} + \ee^2 y_{41} )\}$ are parallel
to the lines  $l_{10}$, $l_{20}$ defined by the pairs $\{(\ee x_{10}
, \ee y_{10} ),(\ee x_{30},
 \ee y_{30}  )\}$ and $\{(\ee x_{20} , \ee y_{20} ),(\ee x_{40},  \ee y_{40}  )\}$. 

Thus the bisector directions between   $l_{11}$, $l_{21}$  in the
plane are the same as for the unperturbed lines  $l_{10}$, $l_{20}$;
so they are obviously parallel to the coordinate axes --- the
principal directions of $\aa(\Omega)$ at $P'$ which is $\ee^2$-close to
both $Z_0= l_{10} \cap l_{20}$ and $Z_1= l_{11} \cap l_{21}$.

Now we have two pairs of circles in $\bbR^3$: $\omega_{10}$,
$\omega_{20}$ passing through the triples of points
 $\{(\ee x_{10} , \ee y_{10} ),P', (\ee x_{30},
 \ee y_{30}  )\}$ and $\{(\ee x_{20} , \ee y_{20} ),P',(\ee x_{40},  \ee y_{40}  )\}$
and the pair $\omega_{11}$, $\omega_{21}$ passing through 
 $\{(\ee x_{10} + \ee^2 x_{11}, \ee y_{10} + \ee^2 y_{11} ),P', (\ee x_{30} + \ee^2 x_{31},
 \ee y_{30} + \ee^2 y_{31}  )\}$ and
 $\{(\ee x_{20} + \ee^2 x_{21}, \ee y_{20} + \ee^2 y_{21} ),P', (\ee x_{40} + \ee^2 x_{41},
 \ee y_{40} + \ee^2 y_{41}  )\}$.  One easily concludes from (\ref{f12}) that  $\omega_{10}$, $\omega_{11}$
are tangent at $P'$, as well as $\omega_{20}$, $\omega_{21}$. Thus
their bisectors are the same and parallel to the axes; the next
$O(\ee^3)$-terms in approximations for the intersection points $A$,
$B$,  $C$,  $D$ will give a $O(\ee^2)$-change in the directions of
$b_1$, $b_2$. $\Box$

\textit{Remark}. From the proof of the previous Theorem we see that
in fact the M\"obius-invariant construction gives the same (with
$O(\ee^2)$-error) principal directions if we will substitute instead
of the original surface $\aa(\Omega)$ its tangent paraboloid at $P$
(after the M\"obius transformation described in the beginning of the
proof). Also it should be noted that the $O(\ee^2)$-approximation
was achieved taking into account approximating paraboloid
of \textit{third order}; the results obtained in this Section are therefore
\textit{third order results} despite the seemingly second-order construction
based on an infinitesimal circle.


\textsl{Proof of Theorem~\ref{Th-pd}.} Now we are in a position to
prove the statement about  $O(\ee^2)$-approximation of principal
directions at $Z_\aa$ in our Euclidean construction. For this we
first perform an auxiliary construction in
$\bbR^3$ (Fig.~\ref{figEu}):

\begin{figure}[htbp]
\begin{center}
\begin{picture}(0,0)%
\includegraphics{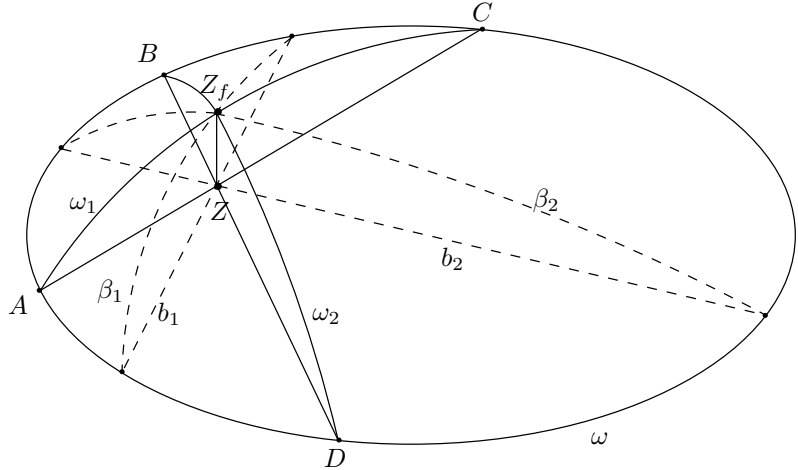}%
\end{picture}%
\setlength{\unitlength}{4144sp}%
\begingroup\makeatletter\ifx\SetFigFont\undefined%
\gdef\SetFigFont#1#2#3#4#5{%
  \reset@font\fontsize{#1}{#2pt}%
  \fontfamily{#3}\fontseries{#4}\fontshape{#5}%
  \selectfont}%
\fi\endgroup%
\begin{picture}(4716,2817)(262,-2631)
\put(262,-1677){\makebox(0,0)[lb]{\smash{{\SetFigFont{10}{12.0}{\familydefault}{\mddefault}{\updefault}{\color[rgb]{0,0,0}$A$}%
}}}}
\put(2151,-2591){\makebox(0,0)[lb]{\smash{{\SetFigFont{10}{12.0}{\familydefault}{\mddefault}{\updefault}{\color[rgb]{0,0,0}$D$}%
}}}}
\put(1468,-1115){\makebox(0,0)[lb]{\smash{{\SetFigFont{10}{12.0}{\familydefault}{\mddefault}{\updefault}{\color[rgb]{0,0,0}$Z$}%
}}}}
\put(1034,-167){\makebox(0,0)[lb]{\smash{{\SetFigFont{10}{12.0}{\familydefault}{\mddefault}{\updefault}{\color[rgb]{0,0,0}$B$}%
}}}}
\put(3038, 78){\makebox(0,0)[lb]{\smash{{\SetFigFont{10}{12.0}{\familydefault}{\mddefault}{\updefault}{\color[rgb]{0,0,0}$C$}%
}}}}
\put(2851,-1388){\makebox(0,0)[lb]{\smash{{\SetFigFont{10}{12.0}{\familydefault}{\mddefault}{\updefault}{\color[rgb]{0,0,0}$b_2$}%
}}}}
\put(631,-1046){\makebox(0,0)[lb]{\smash{{\SetFigFont{10}{12.0}{\familydefault}{\mddefault}{\updefault}{\color[rgb]{0,0,0}$\omega_1$}%
}}}}
\put(2082,-1712){\makebox(0,0)[lb]{\smash{{\SetFigFont{10}{12.0}{\familydefault}{\mddefault}{\updefault}{\color[rgb]{0,0,0}$\omega_2$}%
}}}}
\put(3751,-2453){\makebox(0,0)[lb]{\smash{{\SetFigFont{10}{12.0}{\familydefault}{\mddefault}{\updefault}{\color[rgb]{0,0,0}$\omega$}%
}}}}
\put(1404,-356){\makebox(0,0)[lb]{\smash{{\SetFigFont{10}{12.0}{\familydefault}{\mddefault}{\updefault}{\color[rgb]{0,0,0}$Z_\aa$}%
}}}}
\put(3401,-1036){\makebox(0,0)[lb]{\smash{{\SetFigFont{10}{12.0}{\familydefault}{\mddefault}{\updefault}{\color[rgb]{0,0,0}$\beta_2$}%
}}}}
\put(1157,-1717){\makebox(0,0)[lb]{\smash{{\SetFigFont{10}{12.0}{\familydefault}{\mddefault}{\updefault}{\color[rgb]{0,0,0}$b_1$}%
}}}}
\put(794,-1581){\makebox(0,0)[lb]{\smash{{\SetFigFont{10}{12.0}{\familydefault}{\mddefault}{\updefault}{\color[rgb]{0,0,0}$\beta_1$}%
}}}}
\end{picture}%
\end{center}
\caption{Auxiliary construction} \label{figEu}
\end{figure}


\noindent take the circles $\omega_1$, $\omega_2$ defined by the
point triples $\{A,Z_\aa,C\}$  and  $\{B,Z_\aa,D\}$, then for
the sphere $\Sigma$ where all these points and circles lie we take
its tangent plane $\pi_\Sigma$ at $Z_\aa$. Next we construct
another two circles $\beta_1$, $\beta_2$ as the sections of the
sphere $\Sigma$ by the ``bisectoral planes'' passing through
$Z_\aa$ and the lines $b_1$, $b_2$ respectively. First we remark
that the constructed plane  $\pi_\Sigma$ is  $\ee^2$-close to the
tangent plane  $\pi_\aa$ of $\aa(\Omega)$ at $Z_\aa$: the circles 
$\omega_1$, $\omega_2$  are approximating the osculating circles of
the planar sections of  $\aa(\Omega)$ by the planes  $(A,Z_\aa,C)$ and 
$(B,Z_\aa,D)$. As one can check by a direct computation for any planar
smooth curve $\gamma$, if one takes the circle $\omega$ passing through 
three points lying on $\ee$-distances from each other on $\gamma$, then
the tangents to $\omega$ and $\gamma$ at any of the three pins
will be $\ee^2$-close. So the tangents to the aforementioned planar sections and the tangents to the
respective circles $\omega_1$, $\omega_2$ at $Z_\aa$ are  $\ee^2$-close. (Here
we use the same terminology of ``$\ee^p$-closeness'' for pairs of
lines or planes, this means that the respective angles are
$O(\ee^p)$.) This shows that our sphere $\Sigma$ is  ``very close''
to the sphere $S$ (which should be tangent to $\aa(\Omega)$ at 
$Z_\aa=P$)  used in the M\"obius-invariant construction. We also
see that our Euclidean construction is nearly reduced to the
M\"obius-invariant construction, with the sphere $S$ substituted by
$\Sigma$: the orthogonal projections (along the line
$l=(ZZ_\aa)\perp\pi_{ABCD}$) of $b_1$, $b_2$ onto  $\pi_\Sigma$
are nothing but the tangents to the circles  $\beta_1$, $\beta_2$ 
which are  $\ee^2$-close to the circles which exactly bisect the
angles between the circles   $\omega_1$, $\omega_2$ at  $Z_\aa$. 
This follows from the fact, that  the angle between $\pi_{ABCD}$ and  $\pi_\aa$ is $O(\ee)$:
as one can easily check by a direct computation, in this situation
the (non-infinitesimal) angle between any two lines on one of the planes
and the angle between their orthogonal projections on the other
plane are $\ee^2$-close.

Let us now estimate the angles between the tangents to the circles
$\beta_i$ (lying in $\pi_\Sigma$) and the principal directions of 
$\aa(\Omega)$ at  $Z_\aa$  (lying in $\pi_\aa$, which is 
$\ee^2$-close to  $\pi_\Sigma$). For this we use the same technology
as in the proof of Theorem~\ref{Th-Mic}: first we perform a M\"obius
transformation leaving  $\Sigma$ and $\omega$ invariant and bringing
the point  $Z_\aa$  to the center $O$ of  $\omega$ on $\Sigma$,
then introduce a Cartesian coordinate system with  $O$ as its origin
and $\pi_\Sigma$ being its $(xy)$-plane. Now the equation
(\ref{eq-urp}) of the surface $\aa(\Omega)$ will be modified by small
linear terms: 
\begin{equation}\label{eq-urpsigma}
z= (K_1 x^2 + K_2 y^2)/2 + (a_{30}x^3 + a_{21}x^2 y + a_{12}x y^2 +
a_{03}y^3)
 + \sigma_1x + \sigma_2y + o(\ee^3),
\end{equation}
with $\sigma_i=O(\ee^2)$. The first terms $x_{i0}$, $y_{i0}$ in the
Taylor expansions of the coordinates of the points  $A$, $B$, $C$,
$D$ will be found from the same system (\ref{xy1}) so we have the
relations (\ref{plusminus}). The next terms  $x_{i1}$, $y_{i1}$ are
found from a modified system of the form (\ref{xy2}):
$$
\left\{
\begin{array}{l}
2x_{i0}x_{i1} + 2y_{i0}y_{i1} = 0,\\
K_1 x_{i0}x_{i1} + K_2 y_{i0}y_{i1} + (a_{30}x_{i0}^3 +
a_{21}x_{i0}^2 y_{i0} + a_{12}x_{i0} y_{i0}^2 + a_{03}y_{i0}^3)
  + \frac{\sigma_1}{\ee^2}x_{i0} + \frac{\sigma_2}{\ee^2}y_{i0} = 0.
\end{array}
\right.
$$
We see that the main conclusion (\ref{f12}) still holds; from this
point we just follow the guidelines of the proof of
Theorem~\ref{Th-Mic} and show that the bisectors of the angles
between the circles  $\omega_1$, $\omega_2$ coincide (up to a
$O(\ee^2)$-error) with the directions of the $x$- and $y$-axes. On
the other hand since we know that $ \sigma_i=O(\ee^2)$, using the
standard differential-geometric formulas for calculation of the
coefficients of the second fundamental form and the principal
directions we conclude that the principal directions of   $\aa(\Omega)$ 
given by the equations (\ref{eq-urpsigma}) are $\ee^2$-close to the
$x$- and $y$-axes at $Z_\aa$ as well (as the directions of the
straight lines in $\bbR^3$). $\Box$

\textit{Remark}. From the proof of this Theorem we see that in fact
one can assume (with $O(\ee^2)$-error)  in the Euclidean
construction that the tangent plane at $Z_\aa$ is the readily
constructible plane $\pi_\Sigma$.

\section{Local results}\label{sec-local}

\subsection{Smooth and discrete conjugate quads}\label{sec31}

Suppose that a smooth surface $\aa(\Omega)$ parametrized locally by some 
curvilinear net of conjugate lines is given: $\aa: \Omega
\longrightarrow \bbR^3$, $\Omega \subset \bbR^2 = \{ (u,v)\}$. Take
some initial point $A=\aa(u_0,v_0)$ and points
$B=\aa(u_0+\epsilon,v_0)$, $C=\aa(u_0,v_0+\epsilon)$ at
$\epsilon$-distance on the two conjugate lines of the net on
$\aa(\Omega)$ and  let $D=\aa(u_0+\epsilon,v_0+\epsilon)$ be the
fourth point of the curvilinear quad on $\aa(\Omega)$.

Using the conjugacy condition
\begin{equation}\label{conj}
    \aa_{uv} = a(u,v)\aa_{u} + b(u,v) \aa_{v}
\end{equation}
and its derivatives, one can easily estimate the distance from this
fourth point to the plane $\pi_{ABC}$ defined by $A$, $B$, $C$; we
formulate this as our next Theorem.

We remind that  smooth conjugate nets are supposed to be
\textit{non-degenerate}, that is the directions of the curvilinear
coordinate lines of such a net are non-asymptotic, so the angle
between the tangent vectors $\aa_{u}$ and $\aa_{v}$ is everywhere different from
zero.

\begin{theorem} \label{th-conj}For an arbitrary smooth non-degenerate conjugate net $\aa(\Omega)$, 
$d(D,\pi_{ABC}) =O(\epsilon^4)$.
\end{theorem}

\textsl{Proof.}
    Using the standard Taylor expansions one gets:
\begin{equation}\label{f1}
\begin{array}{l}
    \overrightarrow{AB} = \ee \aa_{u} + \frac{\ee^2}{2}\aa_{uu}
          + \frac{\ee^3}{3!}\aa_{uuu} + \vrln{\vrln{o}}(\ee^3), \\[0.5em]
        \overrightarrow{AC} = \ee \aa_{v} + \frac{\ee^2}{2}\aa_{vv}
          + \frac{\ee^3}{3!}\aa_{vvv} + \vrln{\vrln{o}}(\ee^3),
          \\[0.5em]
        \overrightarrow{AD} = \overrightarrow{AB} + \overrightarrow{AC} + \ee^2\aa_{uv}
          + \frac{\ee^3}{2}(\aa_{uvv} +\aa_{uuv}) + \vrln{\vrln{o}}(\ee^3),
\end{array}
\end{equation}
where all derivatives are taken at the point $(u_0,v_0)$.

Taking into account (\ref{conj}) and its derivatives
$$
\begin{array}{l}
    \aa_{uuv} = a_u\aa_{u} + b_u \aa_{v}
    +  a\aa_{uu}+  b(a\aa_{u}+ b\aa_{v}),
\\[0.5em]
    \aa_{uuv} = a_v\aa_{u} + b_v \aa_{v}
    +  a(a\aa_{u}+ b\aa_{v}) +  b\aa_{vv},
\end{array}
$$
we can compute the triple product $((\overrightarrow{AB}\times
\overrightarrow{AC})\cdot\overrightarrow{AD})$ which gives the
volume of the parallelepiped spanned by the vectors
$\overrightarrow{AB}$, $\overrightarrow{AC}$, $\overrightarrow{AD}$:
$$
  ((\overrightarrow{AB}\times \overrightarrow{AC})\cdot\overrightarrow{AD}) =
  \ee^2 (\overrightarrow{AB}\times \overrightarrow{AC})\cdot \aa_{uv} + \frac{\ee^3}{2}(\overrightarrow{AB}\times \overrightarrow{AC})\cdot(\aa_{uvv}  +\aa_{uuv})
+ \vrln{\vrln{o}}(\ee^5)=
$$
$$
   \frac{\ee^5}{2}\Big[  (\aa_{u} \times \aa_{v})\cdot
        \Big(\aa_{uvv} +\aa_{uuv}\Big) +
       (\aa_{uu} \times \aa_{v})\cdot \aa_{uv} +
       (\aa_{u} \times \aa_{vv})\cdot \aa_{uv}
      \Big] + \vrln{\vrln{o}}(\ee^5)
      =
$$
$$
   \frac{\ee^5}{2}\Big[
   a (\aa_{u} \times \aa_{v})\cdot \aa_{uu}
   + b (\aa_{u} \times \aa_{v})\cdot \aa_{vv}
   + a (\aa_{uu} \times \aa_{v})\cdot \aa_{u}
   + b (\aa_{u} \times \aa_{vv})\cdot \aa_{v}
   \Big] + \vrln{\vrln{o}}(\ee^5)
   =\vrln{\vrln{o}}(\ee^5).
$$
Since the area of the base of this parallelepiped
$|\overrightarrow{AB}\times \overrightarrow{AC}| \simeq \ee^2$, we
get that its height $d = \vrln{\vrln{o}}(\ee^3)$, so at least we may
state that $d =O(\ee^4)$. In fact, using Taylor expansions in
(\ref{f1}) up to order $O(\ee^4)$ one can obtain after a lengthy
computation that
\begin{equation}\label{ee6}
  ((\overrightarrow{AB}\times \overrightarrow{AC})\cdot\overrightarrow{AD}) =
\frac{\ee^6}{12}\Big[(a_u - ab)(\aa_{uu} \times
\aa_{u})\cdot \aa_{v}
  - (b_v - ab) (\aa_{vv} \times \aa_{u})\cdot \aa_{v}
\Big] + \vrln{\vrln{o}}(\ee^6)
\end{equation}
which proves that for a generic conjugate net we have $d \simeq
\ee^4$. $\Box$




In fact we can even choose $M$ (the fourth point of an elementary
planar quad $ABCM$) on the plane $\pi_{ABC}$ sufficiently close to
$D$ and \textit{lying on the given smooth surface} $\aa(\Omega)$: 
\begin{theorem}\label{th-conj-ep3}
For a given smooth non-degenerate conjugate net $\aa: \Omega
\longrightarrow \bbR$ and the plane $\pi_{ABC}$ constructed as above
one can choose a point $M \in \pi_{ABC}\cap \aa(\Omega)$ such that 
$d(M,D) =O(\epsilon^3)$.
\end{theorem}
\textsl{Proof.} As we prove in the Appendix, the intersection of the plane  $\pi_{ABC}$ and the
surface $\aa(\Omega)$ in the $\epsilon$-neighborhood of the points $A$, 
$B$, $C$, $D$ for sufficiently small $\epsilon$ is a curve $I_{ABC}$
close to a quadric -- the Dupin indicatrix of $\aa(\Omega)$,
and  the angle between $\pi_{ABC}$ and $\aa(\Omega)$ is
$\simeq \ee$ in all points of $I_{ABC}$.
According to Theorem~\ref{th-conj},  $d(D,\pi_{ABC}) =O(\epsilon^4)$. From 
this we can easily see that the distance  between $D$ and the
closest to $D$ point $M$ on $I_{ABC}$ should be $O(\ee^4/\ee) =
O(\ee^3)$. $\Box$

\subsection{Curvilinear quads of curvature lines and plane
circular quads}\label{sec32}

For a given smooth surface $\aa(\Omega)$ parametrized by curvature lines 
one can prove a similar result, if we take the circle $\omega_{ABC}$
passing through the points $A$, $B$, $C$ defined as in
Section~\ref{sec31}.

\begin{theorem} \label{th3}
For arbitrary smooth net of curvature lines in a neighborhood of a
non-umbilic point $A$ on $\aa(\Omega)$, the distance between the
fourth point $M$ of intersection of the circle $\omega_{ABC}$ with
$\aa(\Omega)$ and the point $D$ has order 3: $d(D,M)
=O(\epsilon^3)$.
\end{theorem}
In order to prove this theorem we need to establish a few auxiliary
Lemmas.
\begin{lemma} \label{lemma1}
 For arbitrary smooth net of curvature lines
in a neighborhood of a non-umbilic point $A$ on $\aa(\Omega)$, 
$d(D,\omega_{ABC}) =O(\epsilon^3)$.
\end{lemma} 
\textsl{Proof.} Following an approach proposed in \cite{BP96}, we
introduce for each of the points $A$, $B$, $C$, $D$  purely
imaginary quaternion $\mathbf{A}$, $\mathbf{B}$, $\mathbf{C}$,
$\mathbf{D}$. The point $A$ will be chosen as the origin, i.e.\
$\mathbf{A}=0$. The basic quaternions $\mathbf{I}$, $\mathbf{J}$,
$\mathbf{K}$ are chosen to be tangent to the curvature line
directions $\aa_u$, $\aa_v$ and the normal vector to the
surface in the initial point $A$ respectively.
 In view of future applications in the theory of triply orthogonal coordinate systems we include the given surface $\aa(\Omega)$ into  such a system $\aa(u,v,w)$, $\aa(\Omega)$ being one of the coordinate surfaces 
 $\aa(u,v,w=w_0)$. This is always possible (see e.g.\ \cite{Dar}): one can take for example the one-parametric family of surfaces parallel to $\aa(\Omega)$  and the other two one-parametric families of developable surfaces defined by the curvature lines of  $\aa(\Omega)$  and the normals to $\aa(\Omega)$. 

 Introducing for this 3-orthogonal system the Lam\'e coefficients
$H_i(u,v,w)=|\partial_i\aa| $, $\partial_1=\partial/ \partial
u$, $\partial_2=\partial/ \partial v$, $\partial_3=\partial/
\partial w$, normalized vectors $\vec V_i=\partial_i\aa / H_i$
and the rotation coefficients $\beta _{ik}(u) = \partial
_{i}H_{k}/H_{i}$, $i\neq k$, $\beta _{ii}(u) = 0$, we have the
following relations (\cite{Dar}):
\begin{equation}\label{beta}
    \begin{array}{l}
 \partial _{i}H_{k} = \beta _{ik} H_{i},
 \\[0.5em]
 \partial _{i}\vec V_{k} = \beta _{ki} \vec V_{i},
 \\[0.5em]
 \partial _{i}\vec V_{i} =  -\sum_{s \neq i}\beta _{si} \vec V_{s},
\\[0.5em]
 \partial _{j} \beta _{ik}= \beta _{ij} \beta _{jk},
      \quad i\neq j \neq k,
 \\[0.5em]
 \partial _{i}\beta _{ik} + \partial _{k}\beta _{ki} + \sum_{s \neq i,k} \beta _{si}\beta _{sk}=0,
    \end{array}
\end{equation}
where $i,j,k \in \lbrace 1,2,3\rbrace$, $i\neq j \neq k$. In the
initial point $A$ we have $\vec V_1(u_0,v_0)=\mathbf{I}$, $\vec
V_2(u_0,v_0)=\mathbf{J}$, $\vec V_3(u_0,v_0)=\mathbf{K}$. The Taylor
expansions for the other points  are now easily obtained after
differentiation of $\aa_u=H_1\vec V_1$, $\aa_v=H_2\vec V_2$
using (\ref{beta}):
\begin{equation}\label{tay2}
    \begin{array}{l}
\mathbf{B} = \ee H_1 \mathbf{I} + \frac{\ee^2}{2} (\partial_1 H_1
       \mathbf{I} - \beta_{21}H_1 \mathbf{J} -
          \beta_{31}H_1 \mathbf{K}) + \vrln{\vrln{o}}(\ee^2),
\\[0.5em]
\mathbf{C} = \ee H_2 \mathbf{I} + \frac{\ee^2}{2} ( - \beta_{12}H_2
\mathbf{J} + \partial_2 H_2  \mathbf{J} -
          \beta_{32}H_2 \mathbf{K}) + \vrln{\vrln{o}}(\ee^2),
\\[0.5em]
\mathbf{D} = \mathbf{B} + \mathbf{C} +
   \ee^2( \beta_{21}H_2 \mathbf{I} -
          \beta_{12}H_1 \mathbf{J}) + \vrln{\vrln{o}}(\ee^2).
    \end{array}
\end{equation}
Calculating now the quaternionic cross-ratio (\cite{BP96})
$\mathbf{Q}=(\mathbf{A} -\mathbf{B})(\mathbf{B} -\mathbf{C})^{-1}
(\mathbf{C} -\mathbf{D})(\mathbf{D} -\mathbf{A})^{-1}$ and taking
its imaginary part, one can see that $\mathrm{Im}\mathbf{Q}=
\vrln{\vrln{o}}(\ee)$. According to \cite{BP96} we conclude that
$\rho = \vrln{\vrln{o}}(\ee^2)$ for the distance $\rho$ from the
point $D$ to the circle $\omega_{ABC}$ of size  $\simeq \ee$ defined
by the points $A$, $B$, $C$. $\Box$

In the following we fix an $\ee$-patch $\aa(\Omega_\ee)$ of size
$\simeq \ee$  on the surface  $\aa(\Omega)$, such that all the
points $A$, $B$, $C$, $D$ belong to it. All subsequent constructions
will be applied only to this $\ee$-patch $\aa(\Omega_\ee)$.

From Lemma~\ref{lemma2} (proved in the Appendix) we know that the intersection
curve  $I_{ABC}=\aa(\Omega) \cap \pi_{ABC}$ is $\ee^2$-close to the Dupin indicatrix
of $\aa(\Omega)$ at some point $M$. In fact we can take the indicatrix at the following
point $P$ instead:
choose the ``center point'' $P$ on
$\aa(\Omega)$ such that $P$ is the intersection of  $\aa(\Omega)$ and the 
normal to the plane $\pi_{ABC}$ passing through the center of the circle
$\omega_{ABC}$.
In the Cartesian coordinate
system (similar to the system used in the proof of Lemma~\ref{lemma2}),
where the osculating paraboloid at $P$ is ${\cal P}=\{z=(K_1x^2+K_2y^2)/2\}$,
  the points $A$, $B$, $C$, $D$ will have coordinates
\begin{equation}\label{coordABCD}
\begin{array}{l}
 A\Big(-\frac{\ee}{2} + \vrln{\vrln{o}}(\ee),-\frac{\ee}{2} + \vrln{\vrln{o}}(\ee), \frac{K_1+K_2}{8}\ee^2 + \vrln{\vrln{o}}(\ee^2)\Big) , \\
 B\Big(\frac{\ee}{2} + \vrln{\vrln{o}}(\ee),-\frac{\ee}{2} + \vrln{\vrln{o}}(\ee), \frac{K_1+K_2}{8}\ee^2 + \vrln{\vrln{o}}(\ee^2)\Big) , \\
 C\Big(-\frac{\ee}{2} + \vrln{\vrln{o}}(\ee),\frac{\ee}{2} + \vrln{\vrln{o}}(\ee), \frac{K_1+K_2}{8}\ee^2 + \vrln{\vrln{o}}(\ee^2)\Big) , \\
 D\Big(\frac{\ee}{2} + \vrln{\vrln{o}}(\ee),\frac{\ee}{2} + \vrln{\vrln{o}}(\ee), \frac{K_1+K_2}{8}\ee^2 + \vrln{\vrln{o}}(\ee^2)\Big) ,
\end{array}
\end{equation}  so for the angle $\theta $ between the plane $\pi_{ABC}$ and the tangent plane $\pi_P$ to $\aa(\Omega_\ee)$
at $P$ we have $\theta=O(\ee^2)$.
Moreover, the following is true:
\begin{lemma} \label{lemma3}
The centers of the Dupin indicatrix $Q_P={\cal P} \cap \pi_{ABC}$  and $\omega_{ABC}$ are $\ee^2$-close. The
angles of intersection of $\omega_{ABC}$ and $I_{ABC}$  are $\simeq
\ee^0=1$.
\end{lemma}
\textsl{Proof.}  From (\ref{coordABCD}) we easily deduce the first
statement. The radius of $\omega_{ABC}$ is $\ee\sqrt{2} + o(\ee)$,
while the axes of $Q_P$ are $\frac{K_1+K_2}{K_1}+ o(\ee)$ and
$\frac{K_1+K_2}{K_2}+ o(\ee)$. Since we assume $|K_2-K_1|\succ \ee$
on $\Omega$ we see that the intersection angles of $Q_P$ and
$\omega_{ABC}$, so also of  $\omega_{ABC}$ and
 $I_{ABC}$, are $\simeq \ee^0=1$.
$\Box$

Now using the same technique as in the proof of Lemma~\ref{lemma2}, we deduce
that the Dupin indicatrix at $P$ and $I_{ABC}$ are $\ee^2$-close.

\textsl{Proof of Theorem~\ref{th3}.} We know (Theorem~\ref{th-conj})
that $D$ is $\ee^4$-close to  $\pi_{ABC}$. Using the fact that the
angle of intersection of  $\pi_{ABC}$ and $\aa(\Omega_\ee)$ is
$\simeq\ee$ (so  $\succ\ee^2$) we see that $D$ is $\ee^3$-close to
$I_{ABC}= \pi_{ABC} \cap \aa(\Omega)$. On the other hand
(Lemma~\ref{lemma1}) $D$ is $\ee^3$-close to  the circle
$\omega_{ABC}$. Since the angles of intersection of
$\omega_{ABC}$ and $I_{ABC}$  are $\simeq \ee^0=1$, we see 
that $D$ and $M=I_{ABC}\cap \omega_{ABC} $ are $\ee^3$-close  as
well. $\Box$
















\section{Global results}\label{sec-global}

\subsection{Conjugate nets}\label{sec41}

\begin{figure}[htbp]
\begin{center}
\begin{picture}(0,0)%
\includegraphics{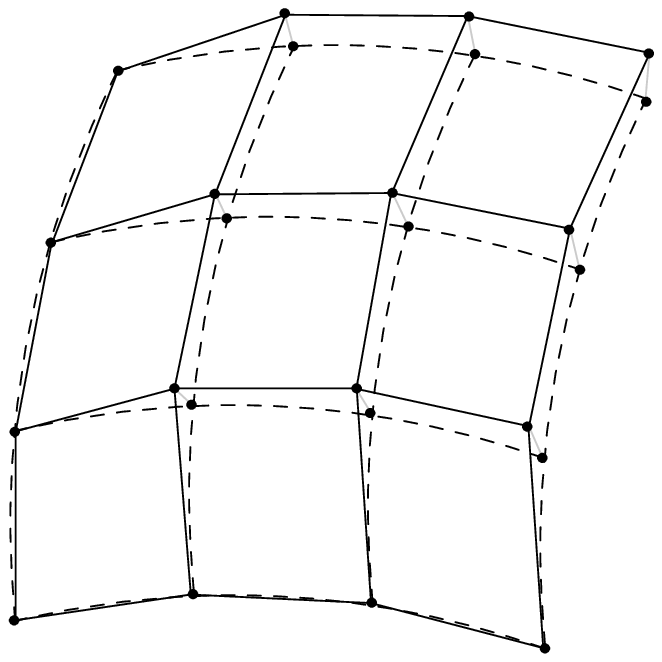}%
\end{picture}%
\setlength{\unitlength}{4144sp}%
\begingroup\makeatletter\ifx\SetFigFont\undefined%
\gdef\SetFigFont#1#2#3#4#5{%
  \reset@font\fontsize{#1}{#2pt}%
  \fontfamily{#3}\fontseries{#4}\fontshape{#5}%
  \selectfont}%
\fi\endgroup%
\begin{picture}(3696,3303)(626,-2525)
\put(956,-2351){\makebox(0,0)[lb]{\smash{{\SetFigFont{9}{10.8}{\familydefault}{\mddefault}{\updefault}{\color[rgb]{0,0,0}$\aa_{00}=\aa^\epsilon_{00}$}%
}}}}
\put(791,-497){\makebox(0,0)[lb]{\smash{{\SetFigFont{9}{10.8}{\familydefault}{\mddefault}{\updefault}{\color[rgb]{0,0,0}$\aa_{02}=\aa^\epsilon_{02}$}%
}}}}
\put(1095,292){\makebox(0,0)[lb]{\smash{{\SetFigFont{9}{10.8}{\familydefault}{\mddefault}{\updefault}{\color[rgb]{0,0,0}$\aa_{03}=\aa^\epsilon_{03}$}%
}}}}
\put(626,-1357){\makebox(0,0)[lb]{\smash{{\SetFigFont{9}{10.8}{\familydefault}{\mddefault}{\updefault}{\color[rgb]{0,0,0}$\aa_{01}=\aa^\epsilon_{01}$}%
}}}}
\put(3487,-2489){\makebox(0,0)[lb]{\smash{{\SetFigFont{9}{10.8}{\familydefault}{\mddefault}{\updefault}{\color[rgb]{0,0,0}$\aa_{30}=\aa^\epsilon_{30}$}%
}}}}
\put(1889,-2259){\makebox(0,0)[lb]{\smash{{\SetFigFont{9}{10.8}{\familydefault}{\mddefault}{\updefault}{\color[rgb]{0,0,0}$\aa_{10}=\aa^\epsilon_{10}$}%
}}}}
\put(2639,-2319){\makebox(0,0)[lb]{\smash{{\SetFigFont{9}{10.8}{\familydefault}{\mddefault}{\updefault}{\color[rgb]{0,0,0}$\aa_{20}=\aa^\epsilon_{20}$}%
}}}}
\put(4322,146){\makebox(0,0)[lb]{\smash{{\SetFigFont{9}{10.8}{\familydefault}{\mddefault}{\updefault}{\color[rgb]{0,0,0}$\aa_{33}$}%
}}}}
\put(2430,682){\makebox(0,0)[lb]{\smash{{\SetFigFont{9}{10.8}{\familydefault}{\mddefault}{\updefault}{\color[rgb]{0,0,0}$\aa^\epsilon_{13}$}%
}}}}
\put(3328,670){\makebox(0,0)[lb]{\smash{{\SetFigFont{9}{10.8}{\familydefault}{\mddefault}{\updefault}{\color[rgb]{0,0,0}$\aa^\epsilon_{23}$}%
}}}}
\put(4205,500){\makebox(0,0)[lb]{\smash{{\SetFigFont{9}{10.8}{\familydefault}{\mddefault}{\updefault}{\color[rgb]{0,0,0}$\aa^\epsilon_{33}$}%
}}}}
\put(3685,-285){\makebox(0,0)[lb]{\smash{{\SetFigFont{9}{10.8}{\familydefault}{\mddefault}{\updefault}{\color[rgb]{0,0,0}$\aa^\epsilon_{32}$}%
}}}}
\put(3473,-1177){\makebox(0,0)[lb]{\smash{{\SetFigFont{9}{10.8}{\familydefault}{\mddefault}{\updefault}{\color[rgb]{0,0,0}$\aa^\epsilon_{31}$}%
}}}}
\put(2374,-498){\makebox(0,0)[lb]{\smash{{\SetFigFont{9}{10.8}{\familydefault}{\mddefault}{\updefault}{\color[rgb]{1,1,1}$\aa_{12}$}%
}}}}
\put(1989,-161){\makebox(0,0)[lb]{\smash{{\SetFigFont{9}{10.8}{\familydefault}{\mddefault}{\updefault}{\color[rgb]{1,1,1}$\aa^\epsilon_{12}$}%
}}}}
\put(2854,-151){\makebox(0,0)[lb]{\smash{{\SetFigFont{9}{10.8}{\familydefault}{\mddefault}{\updefault}{\color[rgb]{1,1,1}$\aa^\epsilon_{22}$}%
}}}}
\put(3181,-558){\makebox(0,0)[lb]{\smash{{\SetFigFont{9}{10.8}{\familydefault}{\mddefault}{\updefault}{\color[rgb]{1,1,1}$\aa_{22}$}%
}}}}
\put(3852,-1492){\makebox(0,0)[lb]{\smash{{\SetFigFont{9}{10.8}{\familydefault}{\mddefault}{\updefault}{\color[rgb]{0,0,0}$\aa_{31}$}%
}}}}
\put(4027,-627){\makebox(0,0)[lb]{\smash{{\SetFigFont{9}{10.8}{\familydefault}{\mddefault}{\updefault}{\color[rgb]{0,0,0}$\aa_{32}$}%
}}}}
\put(2665,283){\makebox(0,0)[lb]{\smash{{\SetFigFont{9}{10.8}{\familydefault}{\mddefault}{\updefault}{\color[rgb]{0,0,0}$\aa_{13}$}%
}}}}
\put(3491,224){\makebox(0,0)[lb]{\smash{{\SetFigFont{9}{10.8}{\familydefault}{\mddefault}{\updefault}{\color[rgb]{0,0,0}$\aa_{23}$}%
}}}}
\put(1847,-1053){\makebox(0,0)[lb]{\smash{{\SetFigFont{9}{10.8}{\familydefault}{\mddefault}{\updefault}{\color[rgb]{0,0,0}$\aa^\epsilon_{11}$}%
}}}}
\put(2244,-1350){\makebox(0,0)[lb]{\smash{{\SetFigFont{9}{10.8}{\familydefault}{\mddefault}{\updefault}{\color[rgb]{0,0,0}$\aa_{11}$}%
}}}}
\put(3041,-1429){\makebox(0,0)[lb]{\smash{{\SetFigFont{9}{10.8}{\familydefault}{\mddefault}{\updefault}{\color[rgb]{0,0,0}$\aa_{21}$}%
}}}}
\put(2674,-1048){\makebox(0,0)[lb]{\smash{{\SetFigFont{9}{10.8}{\familydefault}{\mddefault}{\updefault}{\color[rgb]{0,0,0}$\aa^\epsilon_{21}$}%
}}}}
\end{picture}%
\end{center}
\caption{A discrete conjugate net  $\aa_{ij}^\epsilon$ and a
smooth conjugate net  $\aa_{ij}$ which are $\ee^2$-close. The
points of $\aa_{ij}^\epsilon$ and  $\aa_{ij}$ on the initial
curves coincide} \label{pic411}
\end{figure}

Using the results of Section~\ref{sec31}  one can try to construct
inductively an approximating discrete conjugate net
$\aa^\epsilon$  for sufficiently small  $\epsilon >0$. At the
first glance the following  simplest strategy may be applied:
starting from the initial point $A\equiv \aa_{00}=\aa(u_0,v_0)$ on a
finite piece of a  smooth surface $\aa: \Omega \longrightarrow
\bbR^3$ parametrized by conjugate coordinate lines, fix two series
of  points $\aa_{i,0}= \aa(u_0+i\eps,v_0)$,
 $\aa_{0,i}=\aa(u_0,v_0+i\eps)$
at $\epsilon$-distances on two curvilinear coordinate lines on
$\aa(\Omega)$ passing through $A$.  We can construct the next
point $\aa_{11}^\epsilon$ as the orthogonal projection of
$\aa_{11}=\aa(u_0+\eps,v_0+\eps)$ onto the plane passing through
$\aa_{00}$, $\aa_{10}$, $\aa_{01}$; then
$\aa_{21}^\epsilon$ as the orthogonal projection of
$\aa_{21}=\aa(u_0+2\eps,v_0+\eps)$ onto the plane passing through
$\aa_{10}$, $\aa_{20}$ and  $\aa_{11}^\epsilon$, etc. This
inductive process is shown on Figure~\ref{pic411} with initial
points $\aa_{ij}=\aa(u_0+i\eps,v_0+j\eps)$ on the surface
$\aa(\Omega)$ with smooth coordinate lines (shown as dash-lines).
The approximating discrete conjugate net consists of the points
$\aa_{ij}^\epsilon$,  $\aa_{i0}^\epsilon\equiv\aa_{i0}$,
 $\aa_{0i}^\epsilon\equiv\aa_{0i}$.

~~{\ }



This ``geometric analogue" of the analytic approximation results of
\cite{BMS,M} should give us, if the results of
 \cite{BMS,M} were directly applicable, an estimate $d(\aa_{ij}^\epsilon,\aa_{ij}) \leq C (i+j)\cdot \eps^4$,
so for the complete simple piece  $\aa(\Omega)$ we would have
$d(\aa_{ij}^\epsilon,\aa_{ij}) =O(\eps^3)$. In fact
$d(\aa_{ij}^\epsilon,\aa_{ij})$ shows quadratic behavior:
$d(\aa_{ij}^\epsilon,\aa_{ij}) \leq C i\cdot j\cdot \eps^4$,
so on  $\aa(\Omega)$ we  have  a much weaker estimate
$d(\aa_{ij}^\epsilon,\aa_{ij}) =O(\eps^2)$. Indeed, a
straightforward calculation shows that for the distances
$d(\aa_{ij}^\epsilon,\aa_{ij})$ of the four vertices of any
elementary discrete quadrangle one has
$d(\aa_{i+1,j+1}^\epsilon,\aa_{i+1,j+1}) \simeq
d(\aa_{i+1,j}^\epsilon,\aa_{i+1,j}) +
d(\aa_{i,j+1}^\epsilon,\aa_{i,j+1}) -
d(\aa_{ij}^\epsilon,\aa_{ij}) + \eps^4 \phi(u_0+i\eps,
v_o+j\eps)$ where $\phi(u,v)=(a_u - ab)(\aa_{uu} \times
\aa_{u})\cdot \aa_{v}
  - (b_v - ab) (\aa_{vv} \times \aa_{u})\cdot \aa_{v}$ is the  coefficient of $\eps^6$ in the right hand side of (\ref{ee6}).
Then for example for special surfaces which have constant
$\phi(u,v)\equiv\phi_0\neq 0$ one easily obtains
$d(\aa_{ij}^\epsilon,\aa_{ij}) \simeq i\cdot j\cdot \eps^4$.

Paradoxically, a better strategy consists in choosing the points
$\aa_{ij}^\epsilon$ of the approximating discrete net \textit{on
the given smooth surface $\aa$}, although Theorem~\ref{th-conj-ep3}
suggests that one has worse local error $\simeq\eps^3$ than the
local error $\simeq\eps^4$ of Theorem~\ref{th-conj}. The
construction we propose below is applicable to  arbitrary simple
pieces  $\aa(\Omega)$ \textit{without parabolic points}, i.e.\
points where along with the diagonal coefficient $(\aa_{uv}
\times \aa_{u})\cdot \aa_{v}$ of the second fundamental form
(cf.~(\ref{conj})) at least one of the other coefficients
$(\aa_{uu} \times \aa_{u})\cdot \aa_{v}$, $(\aa_{vv}
\times \aa_{u})\cdot \aa_{v}$ also vanishes.

\begin{figure}[htbp]
\begin{center}
\begin{picture}(0,0)
\includegraphics{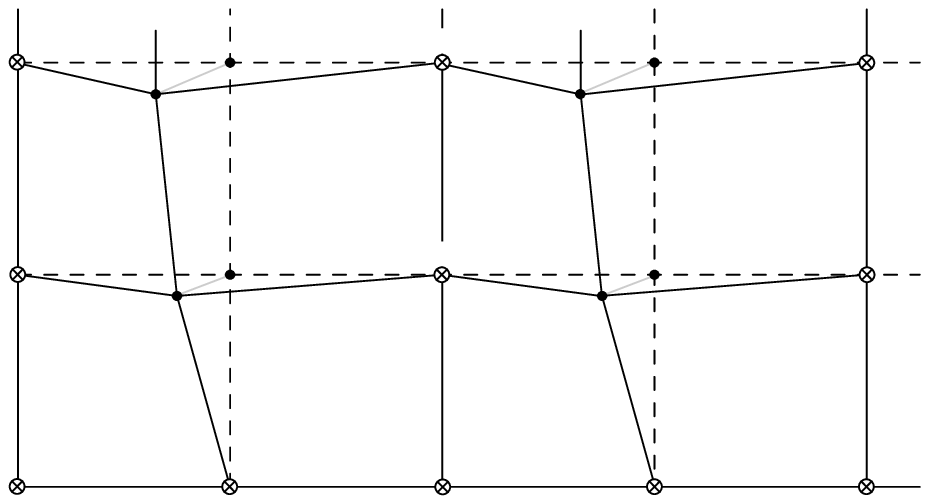}%
\end{picture}%
\setlength{\unitlength}{4144sp}%
\begingroup\makeatletter\ifx\SetFigFont\undefined%
\gdef\SetFigFont#1#2#3#4#5{%
  \reset@font\fontsize{#1}{#2pt}%
  \fontfamily{#3}\fontseries{#4}\fontshape{#5}%
  \selectfont}%
\fi\endgroup%
\begin{picture}(4885,2415)(-112,-1925)
\put(3609,298){\makebox(0,0)[lb]{\smash{{\SetFigFont{9}{10.8}{\familydefault}{\mddefault}{\updefault}{\color[rgb]{0,0,0}$\aa_{32}$}%
}}}}
\put(3611,-673){\makebox(0,0)[lb]{\smash{{\SetFigFont{9}{10.8}{\familydefault}{\mddefault}{\updefault}{\color[rgb]{0,0,0}$\aa_{31}$}%
}}}}
\put(999,-51){\makebox(0,0)[lb]{\smash{{\SetFigFont{9}{10.8}{\familydefault}{\mddefault}{\updefault}{\color[rgb]{0,0,0}$\aa^\epsilon_{12}$}%
}}}}
\put(1118,-984){\makebox(0,0)[lb]{\smash{{\SetFigFont{9}{10.8}{\familydefault}{\mddefault}{\updefault}{\color[rgb]{0,0,0}$\aa^\epsilon_{11}$}%
}}}}
\put(3064,-983){\makebox(0,0)[lb]{\smash{{\SetFigFont{9}{10.8}{\familydefault}{\mddefault}{\updefault}{\color[rgb]{0,0,0}$\aa^\epsilon_{31}$}%
}}}}
\put(2947,-56){\makebox(0,0)[lb]{\smash{{\SetFigFont{9}{10.8}{\familydefault}{\mddefault}{\updefault}{\color[rgb]{0,0,0}$\aa^\epsilon_{32}$}%
}}}}
\put(1668,297){\makebox(0,0)[lb]{\smash{{\SetFigFont{9}{10.8}{\familydefault}{\mddefault}{\updefault}{\color[rgb]{0,0,0}$\aa_{12}$}%
}}}}
\put(1667,-677){\makebox(0,0)[lb]{\smash{{\SetFigFont{9}{10.8}{\familydefault}{\mddefault}{\updefault}{\color[rgb]{0,0,0}$\aa_{11}$}%
}}}}
\put(2269,295){\makebox(0,0)[lb]{\smash{{\SetFigFont{9}{10.8}{\familydefault}{\mddefault}{\updefault}{\color[rgb]{0,0,0}$\aa_{22}=\aa^\epsilon_{22}$}%
}}}}
\put(2269,-680){\makebox(0,0)[lb]{\smash{{\SetFigFont{9}{10.8}{\familydefault}{\mddefault}{\updefault}{\color[rgb]{0,0,0}$\aa_{21}=\aa^\epsilon_{21}$}%
}}}}
\put(316,-1861){\makebox(0,0)[lb]{\smash{{\SetFigFont{9}{10.8}{\familydefault}{\mddefault}{\updefault}{\color[rgb]{0,0,0}$\aa_{00}=\aa^\epsilon_{00}$}%
}}}}
\put(1261,-1861){\makebox(0,0)[lb]{\smash{{\SetFigFont{9}{10.8}{\familydefault}{\mddefault}{\updefault}{\color[rgb]{0,0,0}$\aa_{10}=\aa^\epsilon_{10}$}%
}}}}
\put(2251,-1861){\makebox(0,0)[lb]{\smash{{\SetFigFont{9}{10.8}{\familydefault}{\mddefault}{\updefault}{\color[rgb]{0,0,0}$\aa_{20}=\aa^\epsilon_{20}$}%
}}}}
\put(3241,-1861){\makebox(0,0)[lb]{\smash{{\SetFigFont{9}{10.8}{\familydefault}{\mddefault}{\updefault}{\color[rgb]{0,0,0}$\aa_{30}=\aa^\epsilon_{30}$}%
}}}}
\put(4186,-1861){\makebox(0,0)[lb]{\smash{{\SetFigFont{9}{10.8}{\familydefault}{\mddefault}{\updefault}{\color[rgb]{0,0,0}$\aa_{40}=\aa^\epsilon_{40}$}%
}}}}
\put(-94,-766){\makebox(0,0)[lb]{\smash{{\SetFigFont{9}{10.8}{\familydefault}{\mddefault}{\updefault}{\color[rgb]{0,0,0}$\aa_{01}=\aa^\epsilon_{01}$}%
}}}}
\put(-97,192){\makebox(0,0)[lb]{\smash{{\SetFigFont{9}{10.8}{\familydefault}{\mddefault}{\updefault}{\color[rgb]{0,0,0}$\aa_{02}=\aa^\epsilon_{02}$}%
}}}}
\end{picture}%
\end{center}
\caption{A discrete conjugate net  $\aa_{ij}^\epsilon$ with its
vertices on $\aa(\Omega)$, which are $\ee^2$-close to the vertices of a 
smooth conjugate net $\aa_{ij}$. The original smooth net
$\aa_{ij}$ is schematically shown as a rectangular plane net}
\label{pic412}
\end{figure}

First we fix the vertices $\aa_{i0}^\eps=\aa_{i0}=
\aa(u_0+i\eps,v_0)$ and all even ``columns"
$\aa_{2i,j}^\eps=\aa_{2i,j}= \aa(u_0+2i\eps,v_0+j\eps)$ in the
new approximating discrete conjugate net (marked on
Figure~\ref{pic412} as $\bigotimes$-points).

Then we proceed as follows. Two planes defined by the triples of
points $\{\aa_{00},\aa_{10},\aa_{01}\}$,
$\{\aa_{10},\aa_{20},\aa_{21}\}$ intersect along a straight
line passing through $\aa_{10}^\epsilon=\aa_{10}$. As we show
below in Lemma~\ref{lemma411}, this straight line intersects
$\aa(\Omega)$ in another point $\aa_{11}^\epsilon$ which is close
to $\aa_{11}$, $d(\aa_{11}^\epsilon,\aa_{11}) =O(\ee^3)$,
whereas the distance between this new point $\aa_{11}^\epsilon$
and the line $v=v_0+\ee$ of the smooth conjugate net is of order
$O(\ee^4)$. Other points $\aa_{31}^\epsilon$,
$\aa_{51}^\epsilon$, \ldots of the first row are constructed in
the same way and  obviously have the same error estimates. The
second row $\aa_{2i+1,2}^\epsilon$, $i=0,1,\ldots$ of the
vertices of the approximating discrete net is constructed
analogously using the fixed points  $\aa_{2k,2}^\epsilon$,
$k=0,1,\ldots$ of the second row on $\aa(\Omega)$ and the
constructed complete set of the  points $\aa_{i1}^\epsilon$ of
the first row (cf.~Figure~\ref{pic412}).

As Lemma~\ref{lemma412} below shows, $\aa_{2i+1,2}^\epsilon$ are
now \textit{$\ee^5$-close} to the smooth line $v=v_0+2\ee$ whereas
the distance between $\aa_{2i+1,2}^\epsilon$  and
$\aa_{2i+1,2}$  along this line (the difference of their
$u$-coordinates) is doubled: $\delta u(\aa_{12}^\epsilon)  =
\delta u(\aa_{11}^\epsilon)  + \mu(u_0,v_0+\ee)\ee^3 + o(\ee^3)=
  2  \mu(u_0,v_0)\ee^3 + o(\ee^3)$,
$\delta v(\aa_{12}^\epsilon)  =  - \delta v(\aa_{11}^\epsilon)
+ \nu(u_0,v_0+\ee)\ee^4 + o(\ee^4) =
   o(\ee^4)$,

Here and afterwards we use the following notations: $\delta
v(\aa_{ij}^\epsilon) \equiv v(\aa_{ij}^\epsilon) -
v(\aa_{ij})$, $\delta u(\aa_{ij}^\epsilon) \equiv
u(\aa_{ij}^\epsilon) - u(\aa_{ij}) $, where $u(P)$ and $v(P)$
denote the respective curvilinear coordinates $u$, $v$ for any point
$P$ on the surface $\aa(\Omega)$.

For the third row we will have $\delta v(\aa_{13}^\epsilon)  =
3\mu(u_0,v_0)\ee^3 + o(\ee^3)$, $\delta u(\aa_{13}^\epsilon)  =
\nu(u_0,v_0)\ee^4 + o(\ee^4)$.

This linear "$u$-drift" of the points  $\aa_{1,k}^\epsilon$ in
contrast to "$v$-oscillation" of the same points is easily
explained: the position of the points
$\aa_{i0}^\epsilon=\aa_{i0} = \aa(u_0+i\ee,v_0)$ on the
initial conjugate line $v=v_0$ may be changed using a
reparametrization $ u \longmapsto \bar u = f(u)$ of the  first
conjugate coordinate $u$; thus the $u$-coordinates of the points
$\aa_{i1}^\epsilon$ of the first row may be considered as some
new choice of $u$-parametrization; all subsequent rows simply shift
accordingly.

\begin{figure}[htbp]
\begin{center}
\begin{picture}(0,0)%
\includegraphics{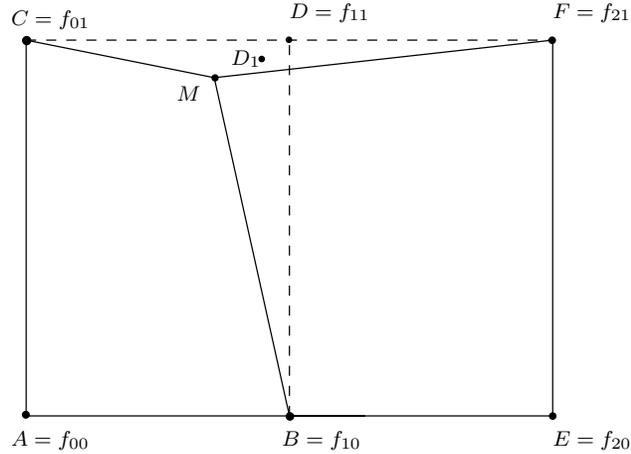}%
\end{picture}%
\setlength{\unitlength}{4144sp}%
\begingroup\makeatletter\ifx\SetFigFont\undefined%
\gdef\SetFigFont#1#2#3#4#5{%
  \reset@font\fontsize{#1}{#2pt}%
  \fontfamily{#3}\fontseries{#4}\fontshape{#5}%
  \selectfont}%
\fi\endgroup%
\begin{picture}(3267,2697)(586,-2077)
\put(586,-2041){\makebox(0,0)[lb]{\smash{{\SetFigFont{8}{9.6}{\familydefault}{\mddefault}{\updefault}{\color[rgb]{0,0,0}$A=\aa_{00}$}%
}}}}
\put(1576, 29){\makebox(0,0)[lb]{\smash{{\SetFigFont{8}{9.6}{\familydefault}{\mddefault}{\updefault}{\color[rgb]{0,0,0}$M$}%
}}}}
\put(586,479){\makebox(0,0)[lb]{\smash{{\SetFigFont{8}{9.6}{\familydefault}{\mddefault}{\updefault}{\color[rgb]{0,0,0}$C=\aa_{01}$}%
}}}}
\put(2251,524){\makebox(0,0)[lb]{\smash{{\SetFigFont{8}{9.6}{\familydefault}{\mddefault}{\updefault}{\color[rgb]{0,0,0}$D=\aa_{11}$}%
}}}}
\put(3826,524){\makebox(0,0)[lb]{\smash{{\SetFigFont{8}{9.6}{\familydefault}{\mddefault}{\updefault}{\color[rgb]{0,0,0}$F=\aa_{21}$}%
}}}}
\put(2206,-2041){\makebox(0,0)[lb]{\smash{{\SetFigFont{8}{9.6}{\familydefault}{\mddefault}{\updefault}{\color[rgb]{0,0,0}$B=\aa_{10}$}%
}}}}
\put(3826,-2041){\makebox(0,0)[lb]{\smash{{\SetFigFont{8}{9.6}{\familydefault}{\mddefault}{\updefault}{\color[rgb]{0,0,0}$E=\aa_{20}$}%
}}}}
\put(1902,243){\makebox(0,0)[lb]{\smash{{\SetFigFont{8}{9.6}{\familydefault}{\mddefault}{\updefault}{\color[rgb]{0,0,0}$D_1$}%
}}}}
\end{picture}%
\end{center}
\caption{The adjacent infinitesimal smooth quads are defined by the
points $A$, $B$, $C$, $D$, $E$ and $F$ (schematically shown as rectangles) and the 
planar quads $ABMC$, $BMFE$} \label{2conjquads}
\end{figure}

\begin{lemma} \label{lemma411}
For two adjacent infinitesimal curvilinear coordinate quadrangles
$ABDC$ and $BDFE$ of $\ee$-size (see Fig.~\ref{2conjquads}) on a
smooth conjugate net $\aa: \Omega \longrightarrow\bbR^3$ without
parabolic points one can find in a $\ee^3$-neighborhood of $D$ a
unique point $M$ on $\aa(\Omega)$ such that the quadrangles $ABMC$
and $BMFE$ are planar. The curvilinear coordinates of $M$ are
$u(M)=u_0+\eps+\delta u_M$, $v(M)=v_0+\eps+\delta v_M$, with $\delta
u_M \sim \mu(u_0,v_0)\ee^3$, $\delta v_M \sim \nu(u_0,v_0)\ee^4$.
\end{lemma}
\textsl{Proof.} We will follow the guidelines of the proof of
Theorem~\ref{th-conj}. Using the Taylor expansions up to order
$O(\ee^4)$ for the Cartesian coordinates of $A(u_0,v_0)$,
$B(u_0+\ee,v_0)$, $C(u_0,v_0+\ee)$, $D(u_0+\ee,v_0+\ee)$,
$E(u_0+2\ee,v_0)$, $F(u_0+2\ee,v_0+\ee)$ and some point
$D_1(u_0+\ee+\delta u_1,v_0+\ee+\delta v_1)$ with some $\delta
u_1=O( \ee^3)$, $\delta v_1=O( \ee^3)$ and taking into consideration
(\ref{conj}), we compute the triple products
$W_1=((\overrightarrow{AB}\times
\overrightarrow{AC})\cdot\overrightarrow{AD_1}) = \ee^6
\gamma(u_0,v_0) + \frac{\ee^3}{2}\Big(\delta u_1  (\aa_{uu}
\times \aa_{u})\cdot \aa_{v}
 + \delta v_1  (\aa_{vv} \times \aa_{u})\cdot \aa_{v} \Big) + o(\ee^6)$,
$W_2=((\overrightarrow{BE}\times
\overrightarrow{BD_1})\cdot\overrightarrow{FB})= \ee^6
\gamma(u_0,v_0) + \frac{\ee^3}{2}\Big(\delta u_1  (\aa_{uu}
\times \aa_{u})\cdot \aa_{v}
 - \delta v_1  (\aa_{vv} \times \aa_{u})\cdot \aa_{v} \Big) + o(\ee^6)$ where $\gamma(u,v)$
is a combination of the coefficients $\aa(u,v)$, $\beta(u,v)$ of
(\ref{conj}), their derivatives and the triple products $
(\aa_{uu} \times \aa_{u})\cdot \aa_{v}$,  $(\aa_{vv}
\times \aa_{u})\cdot \aa_{v}$. So if both of the  latter
triple products do not vanish, one can choose unique $\delta u_1 =
\bar \mu(u_0,v_0)\ee^3$, $\delta v_1 =0$ such that $W_1 =O( \ee^7)$,
$W_2 =O( \ee^7)$.

This shows (as in the proof of Theorem~\ref{th-conj-ep3})  that
$D_1$ lies in $\ee^4$- neighborhood of the curves $I_{ABC}$,
$I_{BEF}$ of intersection of $\aa(\Omega)$ with the planes $(ABC)$,
$(BEF)$. Since the angle of intersection of $I_{ABC}$, $I_{BEF}$ is
$\simeq 1$, we conclude that the point $M$ of their intersection,
close to $D_1$, is in fact $\ee^4$-close to $D_1$.

Using Taylor expansions for the same points up to order $O(\ee^5)$,
one can prove that for a generic surface actually $\delta u_M \sim
\mu(u_0,v_0)\ee^3$, $\delta v_M \sim \nu(u_0,v_0)\ee^4$ with $\mu
\not \equiv 0$, $\nu \not \equiv 0$. $\Box$

If one repeats  the same calculation using instead of $B$ a point
$B_1(u_0+\ee+\delta u(B_1), v_0+\delta v(B_1))$ with some $\delta
u(B_1) =O( \ee^2)$, $\delta v(B_1) =O( \ee^3)$, one will get from
the conditions $W_1 =O( \ee^7)$, $W_2 =O( \ee^7)$ that $\delta u_1 =
\delta u(B_1) + \mu(u_0,v_0)\ee^3 + o(\ee^3)$, $\delta v_1 = -
\delta v(B_1) +  \nu(u_0,v_0)\ee^4
 +  \lambda(u_0,v_0) \delta v(B_1)\ee +  \gamma(u_0,v_0) \delta u(B_1)\ee^2 + o(\ee^4)$
where $\mu(u_0,v_0)$, $\nu(u_0,v_0)$, $\lambda(u_0,v_0)$, $ \gamma
(u_0,v_0)$ are some algebraic combinations of the coefficients
$\aa(u,v)$, $\beta(u,v)$ of (\ref{conj}), their derivatives and the
triple products $ (\aa_{uu} \times \aa_{u})\cdot
\aa_{v}$,  $(\aa_{vv} \times \aa_{u})\cdot \aa_{v}$.
Hence the following estimates hold for the second point $M$ of
intersection of $I_{AB_1C}$ and $I_{B_1EF}$:
\begin{lemma}  \label{lemma412}
For two adjacent infinitesimal curvilinear coordinate quadrangles
$ABDC$ and $BDFE$ of $\ee$-size on a smooth conjugate net $\aa:
\Omega \longrightarrow\bbR^3$ without parabolic points and some
point $B_1(u_0+\ee+\delta u(B_1), v_0+\delta v(B_1))$,  $\delta
u(B_1) =O( \ee^2)$, $\delta v(B_1) =O( \ee^3)$, one can find in a
$\ee^2$-neighborhood of $D$ a unique point $M$ on $\aa(\Omega)$ such
that the quadrangles $AB_1MC$ and $B_1MFE$ are planar. The
curvilinear coordinates of $M$ are $u(M)=u_0+\eps+\delta u_M$,
$v(M)=v_0+\eps+\delta v_M$, with $\delta u_M = \delta u(B_1) +
\mu(u_0,v_0)\ee^3 + o(\ee^3)$, $\delta v_M = - \delta v(B_1) +
\nu(u_0,v_0)\ee^4  +  \lambda(u_0,v_0)\delta v(B_1)\ee  +
\gamma(u_0,v_0) \delta u(B_1)\ee^2  + o(\ee^4)$.
\end{lemma}
\textit{Remark}. As before, one can explain the
``$\ee^2$-tolerance'' along the $u$-lines simply by removing this
$u$-shift using a reparametrization of the curvilinear coordinate
$u$. Note that we need this reparametrization only locally, for a
\textit{given} pair of elementary infinitesimal quadrangles, in
order to establish this estimate for them, and not globally, for the
complete discrete conjugate net.

For global estimates we need the following simple result
(\cite{BMS}):
\begin{lemma}\label{lemmaGronwall}(Discrete Gr\"onwall estimate) Assume that a nonnegative function $\Delta: \bbN
 \longrightarrow \bbR$ satisfies
\begin{equation}\label{Gronwall0}
\Delta(n+1) \leq (1+\eps K)\Delta(n) + \kappa
\end{equation}
with nonnegative constants $K$, $\kappa$. Then
\begin{equation}\label{Gronwall}
\Delta(n) \leq (\Delta(0) + n\kappa)\exp(Kn\ee).
\end{equation}
\end{lemma}
\begin{theorem} For a
smooth conjugate net without parabolic points $\aa:\Omega\to
\bbR^3$
and sufficiently small $\epsilon >0$, there exists a discrete
conjugate net $\aa^\epsilon$ with all its points
$\aa_{ij}^\epsilon$
 on $\aa(\Omega)$, such that $d(\aa_{ij},\aa_{ij}^\epsilon)
=O( \epsilon^2)$.
\end{theorem}
\textsl{Proof.} From Lemma~\ref{lemma412} we immediately see that
for each column $k=2i+1$ the function $\Delta(n)=|\delta
u(\aa_{kn}^\epsilon)|$ satisfies the estimate (\ref{Gronwall0})
with $K=0$, $\kappa= 2 \ee^3\max_{(u,v) \in \Omega} |\mu(u,v)| $ for
sufficiently small $\ee$, so for all $n$ we have  (\ref{Gronwall}).
Since the number of steps in each column is $\simeq 1/\ee$ we obtain
the global estimate  $|\delta u(\aa_{2i+1,j}^\epsilon)| \leq
C\ee^2$.

Now we set  $\Delta(n)=|\delta v(\aa_{kn}^\epsilon)|$ for each
column $k=2i+1$. Taking the estimates of Lemma~\ref{lemma412} and
the already obtained global estimate for   $|\delta
u(\aa_{2i+1,j}^\epsilon)|$ we get (\ref{Gronwall0}) with
$K=\max_{(u,v) \in \Omega} |\lambda(u,v)| $, $\kappa= 2
\ee^4\max_{(u,v) \in \Omega}
  ( |\mu(u,v)| +  |\gamma(u,v)|) $, so  (\ref{Gronwall}) gives the global estimate
$|\delta v(\aa_{2i+1,j}^\epsilon)| \leq C\ee^3$.
$\Box$

\textit{Remark}. As we have observed above, the $v$-shift of the
points $\aa_{2i+1,j}^\epsilon$ has a much more interesting
oscillating behavior; in particular one may conjecture  that
$|\delta(\aa_{2i+1,j}^\epsilon)| =O( \ee^4)$ globally.

\subsection{Curvature lines and circular nets}\label{sec42}

\begin{figure}[htbp]
\begin{center}
\begin{picture}(0,0)%
\includegraphics{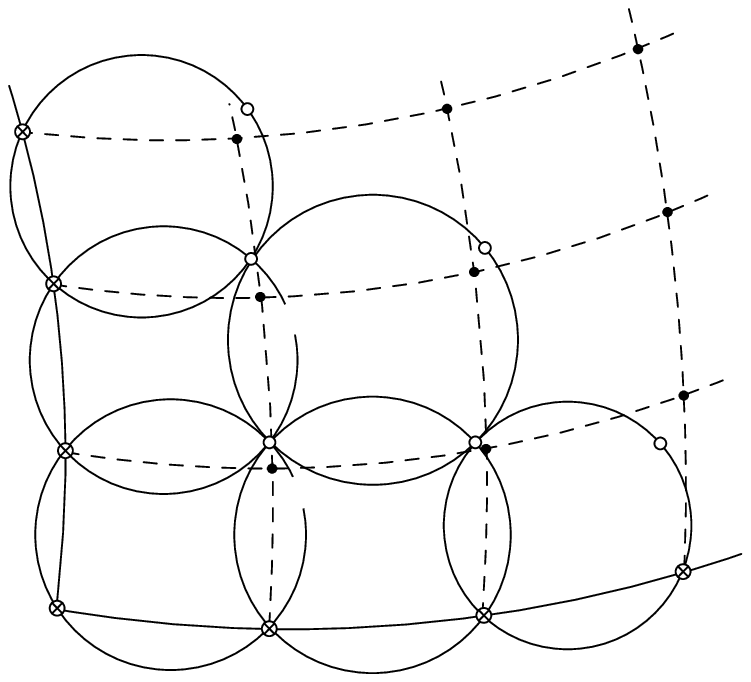}%
\end{picture}%
\setlength{\unitlength}{4144sp}%
\begingroup\makeatletter\ifx\SetFigFont\undefined%
\gdef\SetFigFont#1#2#3#4#5{%
  \reset@font\fontsize{#1}{#2pt}%
  \fontfamily{#3}\fontseries{#4}\fontshape{#5}%
  \selectfont}%
\fi\endgroup%
\begin{picture}(4040,3177)(339,-5817)
\put(4162,-4516){\makebox(0,0)[lb]{\smash{{\SetFigFont{10}{12.0}{\familydefault}{\mddefault}{\updefault}{\color[rgb]{1,1,1}$\aa_{31}$}%
}}}}
\put(3773,-4776){\makebox(0,0)[lb]{\smash{{\SetFigFont{10}{12.0}{\familydefault}{\mddefault}{\updefault}{\color[rgb]{1,1,1}$\aa^\epsilon_{31}$}%
}}}}
\put(2194,-3060){\makebox(0,0)[lb]{\smash{{\SetFigFont{10}{12.0}{\familydefault}{\mddefault}{\updefault}{\color[rgb]{1,1,1}$\aa^\epsilon_{13}$}%
}}}}
\put(1811,-3388){\makebox(0,0)[lb]{\smash{{\SetFigFont{10}{12.0}{\familydefault}{\mddefault}{\updefault}{\color[rgb]{1,1,1}$\aa_{13}$}%
}}}}
\put(3237,-3685){\makebox(0,0)[lb]{\smash{{\SetFigFont{10}{12.0}{\familydefault}{\mddefault}{\updefault}{\color[rgb]{0,0,0}$\aa^\epsilon_{22}$}%
}}}}
\put(2902,-4032){\makebox(0,0)[lb]{\smash{{\SetFigFont{10}{12.0}{\familydefault}{\mddefault}{\updefault}{\color[rgb]{0,0,0}$\aa_{22}$}%
}}}}
\put(2225,-4093){\makebox(0,0)[lb]{\smash{{\SetFigFont{9}{10.8}{\familydefault}{\mddefault}{\updefault}{\color[rgb]{0,0,0}$\aa_{12}$}%
}}}}
\put(2275,-4885){\makebox(0,0)[lb]{\smash{{\SetFigFont{9}{10.8}{\familydefault}{\mddefault}{\updefault}{\color[rgb]{0,0,0}$\aa_{11}$}%
}}}}
\put(1811,-3844){\makebox(0,0)[lb]{\smash{{\SetFigFont{9}{10.8}{\familydefault}{\mddefault}{\updefault}{\color[rgb]{0,0,0}$\aa^\epsilon_{12}$}%
}}}}
\put(3286,-4757){\makebox(0,0)[lb]{\smash{{\SetFigFont{10}{12.0}{\familydefault}{\mddefault}{\updefault}{\color[rgb]{0,0,0}$\aa_{21}$}%
}}}}
\put(2854,-4632){\makebox(0,0)[lb]{\smash{{\SetFigFont{9}{10.8}{\familydefault}{\mddefault}{\updefault}{\color[rgb]{0,0,0}$\aa^\epsilon_{21}$}%
}}}}
\put(2957,-5725){\makebox(0,0)[lb]{\smash{{\SetFigFont{9}{10.8}{\familydefault}{\mddefault}{\updefault}{\color[rgb]{0,0,0}$\aa_{20}=\aa^\epsilon_{20}$}%
}}}}
\put(613,-5565){\makebox(0,0)[lb]{\smash{{\SetFigFont{9}{10.8}{\familydefault}{\mddefault}{\updefault}{\color[rgb]{0,0,0}$\aa_{00}=\aa^\epsilon_{00}$}%
}}}}
\put(4157,-5352){\makebox(0,0)[lb]{\smash{{\SetFigFont{9}{10.8}{\familydefault}{\mddefault}{\updefault}{\color[rgb]{0,0,0}$\aa_{30}=\aa^\epsilon_{30}$}%
}}}}
\put(1886,-5781){\makebox(0,0)[lb]{\smash{{\SetFigFont{9}{10.8}{\familydefault}{\mddefault}{\updefault}{\color[rgb]{0,0,0}$\aa_{10}=\aa^\epsilon_{10}$}%
}}}}
\put(531,-4708){\makebox(0,0)[lb]{\smash{{\SetFigFont{9}{10.8}{\familydefault}{\mddefault}{\updefault}{\color[rgb]{0,0,0}$\aa_{01}=\aa^\epsilon_{01}$}%
}}}}
\put(478,-3933){\makebox(0,0)[lb]{\smash{{\SetFigFont{9}{10.8}{\familydefault}{\mddefault}{\updefault}{\color[rgb]{0,0,0}$\aa_{02}=\aa^\epsilon_{02}$}%
}}}}
\put(339,-3248){\makebox(0,0)[lb]{\smash{{\SetFigFont{9}{10.8}{\familydefault}{\mddefault}{\updefault}{\color[rgb]{0,0,0}$\aa_{03}=\aa^\epsilon_{03}$}%
}}}}
\put(2298,-4654){\makebox(0,0)[lb]{\smash{{\SetFigFont{9}{10.8}{\familydefault}{\mddefault}{\updefault}{\color[rgb]{0,0,0}$\aa^\epsilon_{11}$}%
}}}}
\end{picture}%
\end{center}
\caption{A circular net  $\aa_{ij}^\epsilon$ with its vertices on
$\aa(\Omega)$, which are $\ee^2$-close to the vertices of the  smooth 
net of conjugate lines $\aa_{ij}$ } \label{pic421}
\end{figure}

We start the construction of a circular discrete approximation of a
given finite piece $\aa(\Omega)$ of a smooth surface without umbilic
points parametrized by curvature lines fixing points
$\aa_{i0}^\epsilon\equiv \aa_{i0}= \aa(u_0+i\eps,v_0)$,
$\aa_{0i}^\epsilon\equiv\aa_{0i}=\aa(u_0,v_0+i\eps)$ on two
initial curvature lines passing throw a point $\aa_{00}=
\aa(u_0,v_0)$ (cf.~Figure~\ref{pic421}).


~~{\ }


\vspace{0.2em}

The first circle is defined by the triple $\{\aa_{00},
\aa_{10}, \aa_{01}\}$, its fourth point of intersection with
the surface will be chosen as the next constructed point
$\aa_{11}^\ee$ of the discrete circular net. On the next step we
can define two circles by the triples of points $\{\aa_{01},
\aa_{02}, \aa_{11}^\eps\}$ and
 $\{\aa_{10}, \aa_{20}, \aa_{11}^\eps\}$, the fourth points of intersection of these circles with the surface will give us
 $\aa_{12}^\ee$ and  $\aa_{21}^\ee$. The behavior of the new points on the third step is of crucial importance.
Their error estimates can be obtained from the following
generalization of Theorem~\ref{th3}:
\begin{lemma}  \label{lemma421}
Let the points $A$, $B$ and $C$ on $\aa(\Omega)$ be chosen, such
that $u(A)=u_0+\delta u_A$, $v(A)=v_0+\delta v_A$,
$u(B)=u_0+\ee+\delta u_B$, $v(B)=v_0+\delta v_B$, $u(C)=u_0+\delta
u_C$, $v(C)=v_0+\ee+\delta v_C$, with all $\delta u$, $\delta v$ of
order $=O(\ee^3)$. Then the curvilinear coordinates of the fourth
point $M$ of intersection of the circle $\omega_{ABC}$ with
$\aa(\Omega)$ are: $u(M)=u(\aa_{11})+\delta u_M$,
$v(M)=v(\aa_{11})+\delta v_M$, with
\begin{equation}\label{f421}
\begin{array}{ll}
\delta u_M = &\delta u_A +\delta u_B -\delta u_C
  +\ee(\delta v_B -\delta v_A)\mu_1(u_0,v_0)\\
 &+\ee(\delta u_C -\delta u_A)\mu_2(u_0,v_0)
 +\mu_3(u_0,v_0)\ee^3 +\mu_4(u_0,v_0)\ee^4 + o(\ee^4),  \\
\delta v_M = &\delta v_A -\delta v_B +\delta v_C
  +\ee(\delta v_B -\delta v_A)\nu_1(u_0,v_0) \\
 &+\ee(\delta u_C -\delta u_A)\nu_2(u_0,v_0)
 +\nu_3(u_0,v_0)\ee^3 +\nu_4(u_0,v_0)\ee^4 + o(\ee^4),
\end{array}
\end{equation}
\end{lemma}
\textsl{Proof.}  We will use the same quaternionic cross-ratio
$\mathbf{Q}=(\mathbf{A} -\mathbf{B})(\mathbf{B} -\mathbf{C})^{-1}
(\mathbf{C} -\mathbf{D})(\mathbf{D} -\mathbf{A})^{-1}$. The obvious
changes are necessary in the Taylor expansions (\ref{tay2}), where
we shall take into consideration the $\delta u$- and $\delta
v$-shifts of the points; the expansions themselves should include
terms up to order $O(\ee^4)$. For simplicity we will assume $\delta
u_A=0$, $\delta v_A=0$, the general case requires only a shift by
$\delta u_A$, $\delta v_A$ in the net of curvature lines on
$\aa(\Omega)$.

Fixing for the moment the points $A$, $B$, $C$ and introducing  $D
\in \aa(\Omega)$, $u(D)=u_0+\ee +\delta u_D$, $v(D)=v_0+\ee+\delta
v_D$ with some $\delta u_D=O(\ee^3)$, $\delta v_D=O(\ee^3)$,
one can compute the imaginary part of $\mathbf{Q}$:
$$
\begin{array}{c}
Im(\mathbf{Q}) =
\mathbf{I} \big(\rho_1(u_0,v_0) \ee^3 +   \theta H_2(\delta v_D   -   \delta v_C   +   \delta v_B)   + o(\ee^3)\big) \\[1em]
   +  \mathbf{J} \big(\rho_2(u_0,v_0)\ee^3 +   \theta H_1( - \delta u_D  -   \delta u_C   +   \delta u_B)  + o(\ee^3) \big) \\[1em]
+  \mathbf{K}  \big( \rho_3\ee^2 + \ee^{-1}H_1H_2
    ( H_1^2( \delta u_B-\delta u_D -   \delta u_C) + H_2^2 (\delta v_D   +   \delta v_C -   \delta v_B )  )
   + o(\ee^2)\big), \\[1em]
\end{array}
$$
where $\theta = \frac{1}{2}(H_1^2H_2\beta_{32} - H_1H_2^2\beta_{31})
= H_1^2H_2^2(K_2-K_1)/2$, $K_i$ being the principal curvatures at
the point $A$. Since we assume $K_1 \neq K_2$, choosing
$$
\begin{array}{l}
\delta u_D  =  \delta u_B   -   \delta u_C + \ee^3\mu_3(u_0,v_0) + o(\ee^3),\\
\delta v_D  = -\delta v_B   +   \delta v_C + \ee^3\nu_3(u_0,v_0) +
o(\ee^3),
\end{array}
$$
one will achieve $Im(\mathbf{Q}) = o(\ee^3) \mathbf{I}  + o(\ee^3)
\mathbf{J} + o(\ee^2) \mathbf{K}$. As a lengthier computation shows,
adding appropriate fourth-order corrections in Taylor expansions, 
given in (\ref{f421}) (set there $\delta u_A=0$, $\delta v_A=0$), we
get $Im(\mathbf{Q}) = o(\ee^4) \mathbf{I}  + o(\ee^4) \mathbf{J} +
o(\ee^3) \mathbf{K}$. This means that the point $D$ with $\delta
u_D$, $\delta v_D$ given in (\ref{f421}) is $\ee^6$-close to the
plane $(ABC)$ (this can be also checked using the triple product
$(\overrightarrow{AB}\times
\overrightarrow{AC})\cdot\overrightarrow{AD}$) and $\ee^5$-close to
the circle $\omega_{ABC}$ passing through $A$, $B$, $C$. Using the
the fact that the angle between $\aa(\Omega)$ and the plane $(ABC)$
is $\simeq \ee$ at $D$, we conclude that for the fourth point $M$ of
intersection of  $\omega_{ABC}$ and $\aa(\Omega)$ we can keep the
same expansions (\ref{f421}). $\Box$

\textit{Remark}. Linear behavior of  $\delta u_D$ and $\delta v_D$
in   (\ref{f421}) is valid precisely in orders $O(\ee^3)$ and
$O(\ee^4)$, corrections of order  $O(\ee^5)$ are nonlinear w.r.t.\
$\delta u$- and $\delta v$-shifts of the points
 $A$, $B$ and $C$.

Using  (\ref{f421}) we calculate $
   \delta u (\aa^\ee_{12}) = \delta u(\aa_{01}) +\delta u(\aa_{11}^\ee) -\delta u(\aa_{02})
  +\mu_3(u_0,v_0+\ee)\ee^3  + o(\ee^3) = 2\mu_3(u_0,v_0)\ee^3  + o(\ee^3)
$, $     \delta v (\aa^\ee_{12}) = \delta v_(\aa_{01}) -\delta
v(\aa_{11}^\ee) +\delta v(\aa_{02})
  +\ee(\delta v(\aa_{11}^\ee) -\delta v(\aa_{01}))\nu_1(u_0,v_0+\ee)
 +\ee(\delta u(\aa_{02}) -\delta u(\aa_{01}))\nu_2(u_0,v_0+\ee)
 +\nu_3(u_0,v_0+\ee)\ee^3 +\nu_4(u_0,v_0+\ee)\ee^4 + o(\ee^4) =
\ee^4(\partial_v\nu_3(u_0,v_0) + \nu_1(u_0,v_0)\nu_3(u_0,v_0)) +
o(\ee^4) $. As we see now, for points $\aa^\ee_{12}$ (respectively
$\aa^\ee_{21}$) only the shifts \textit{along} the respective
curvature lines is of order $O(\ee^3)$ (only for surfaces of special
type this shift may degenerate to $0$), while the
\textit{perpendicular} shift is of order $O(\ee^4)$. For
$\aa^\ee_{22}$ we now obtain a remarkable result: \textit{both its
shifts are of order $O(\ee^4)$}, compared to the error estimates of
order $O(\ee^3)$ for  $\aa^\ee_{11}$: $     \delta u (\aa^\ee_{22})
= 2\ee^4(\partial_u\mu_3(u_0,v_0) + \mu_2(u_0,v_0)\mu_3(u_0,v_0)) +
o(\ee^4) $, $     \delta v (\aa^\ee_{22}) =
2\ee^4(\partial_v\nu_3(u_0,v_0) + \nu_1(u_0,v_0)\nu_3(u_0,v_0)) +
o(\ee^4) $.

This observation suggests us to partition the complete discrete
lattice $\aa^\ee_{ij}$ into 3 sublattices:

a) the even sublattice of points $\aa^\ee_{2i,2j}$,

b) the odd sublattice of points $\aa^\ee_{2i+1,2j+1}$,

c) the intermediate sublattice of points $\aa^\ee_{2i+1,2j}$,
$\aa^\ee_{2i,2j+1}$.

Infinitesimally, in a $(N\ee)$-neighborhood of the initial point
$\aa^\ee_{00}$  (with $N \ll 1/\ee$), the even sublattice has both
$\delta u_D$- and $\delta v_D$-shifts of order $O(\ee^4)$, for the
odd sublattice they have order $O(\ee^3)$, for the intermediate one
only the shifts along the respective curvature lines are of order
$O(\ee^3)$, the perpendicular shifts being of order $O(\ee^4)$
again.




Accumulation of errors for these sublattices is also different:
\textit{linear} accumulation of $O(\ee^3)$-errors for the odd
sublattice, \textit{quadratic} for the even sublattice and mixed for
the intermediate sublattice: {linear} accumulation of
$O(\ee^3)$-shifts along the respective curvature lines and quadratic
in perpendicular direction: $   \delta u (\aa^\ee_{2i,2j}) \simeq
i\cdot j\cdot \eps^4$, $   \delta v (\aa^\ee_{2i,2j}) \simeq
i\cdot j\cdot \eps^4$, $   \delta u (\aa^\ee_{2i+1,2j+1}) \simeq  (
i+ j)\cdot \eps^3$, $   \delta v (\aa^\ee_{2i+1,2j+1}) \simeq  ( i+
j)\cdot \eps^3$, $\delta u (\aa^\ee_{2i,2j+1}) \simeq   i\cdot
j\cdot \eps^4$, $   \delta v (\aa^\ee_{2i,2j+1}) \simeq   (i+
j)\cdot \eps^3$.

These estimates are easily  obtained using (\ref{f421}) (see below
the proof of Theorem~\ref{th4222}), for example the first one
follows from:
\begin{equation}\label{evenlat}
\begin{array}{l}
    \delta u (\aa^\ee_{2(i+1),2(j+1)}) =  \delta u (\aa^\ee_{2(i+1),2j})
   +  \delta u (\aa^\ee_{2i,2(j+1)}) -  \delta u (\aa^\ee_{2i,2j}) + {}\\[0.5em]
{~} \quad 2\ee^4(\partial_u\mu_3 + \mu_2\mu_3) - 2\ee\mu_2(\delta u
(\aa^\ee_{2i,2(j+1)}) - \delta u (\aa^\ee_{2i,2j})) + o(\ee^4).
\end{array}\end{equation}
Note that the estimates of Lemma~\ref{lemma421} do not require the
initial three points $A$, $B$, $C$ to lie near the grid points
$\aa_{ij}$ of the original smooth net and the distances $d(A,B)$,
$d(A,C)$ need not to be equal, both just have to
be of order $O(\ee)$. 
\begin{theorem}\label{th4222}
For a smooth surface without umbilic points parametrized by
curvature lines  $\aa:\Omega\to \bbR^3$, $\Omega= \{(i,v)|
u^2+v^2<1\}$   and sufficiently small $\epsilon >0$, there exists a
discrete circular net $\aa^\epsilon$ with all its points
$\aa_{ij}^\epsilon$ on $\aa(\Omega)$, such that
$d(\aa_{ij},\aa_{ij}^\epsilon) =O( \epsilon^2)$.
\end{theorem}
\textsl{Proof.} We will give the details for the global estimate $
\delta u (\aa^\ee_{2i,2j}) \leq  C\cdot  i\cdot j\cdot \eps^4$, the
other are proved in the same way.

First define a function $S(k,n)=|\delta
u(\aa_{2k,2(n+1)}^\epsilon) - \delta
u(\aa_{2k,2n}^\epsilon)|$. From (\ref{evenlat}) we see that
$$ S(k+1,n) \leq S(k,n) + 2\ee|\mu_2|S(k,n) + 2\ee^4|\partial_u\mu_3 + \mu_2\mu_3| +o(\ee^4)
$$
and using Lemma~\ref{lemmaGronwall} we immediately obtain
$$S(k,n) \leq C \ee^4 k.
$$
Now for the functions $\Delta_k(n)=|\delta
u(\aa_{2k,2n}^\epsilon)|$ taken separately for each column of the
even sublattice one gets from  (\ref{evenlat}) that
$$\Delta_k(n+1) \leq \Delta_k(n) + (1+ 2\ee|\mu_2|)S(k,n) + 2\ee^4 K \leq
  \Delta_k(n) + \ee^4 (2Ck +K)
$$
with $K=2\max_{(u,v) \in \Omega} |\partial_u\mu_3 + \mu_2\mu_3| $.
So  (\ref{Gronwall}) gives us the global estimate $ \Delta_k(n) \leq
\bar C\cdot  k\cdot n\cdot \eps^4$ $\Box$

In \cite{BMS} one can find a similar result with the same order
$\epsilon^2$ of approximation but without the requirement
$\aa_{ij}^\epsilon \in \aa(\Omega)$.

\section*{Acknowledgments}
The authors wish to thank Dr. I.~Dynnikov for the idea of simplification 
of the proof of Lemma~\ref{thd} and Prof.~P.~Schr\"oder for valuable discussions.

\appendix

\section*{Appendix}\label{app}  

\begin{lemma} \label{lemma2} Let an $\ee$-neighborhood 
$\aa: \Omega_\ee \rightarrow \bbR^3$ of a point $M$ on a smooth
surface be given, and for the principal curvatures one has
$K_1^2+K_2^2 >0$ at $M$. Then the intersection
$I_{ABC}$ of a plane $\pi_{ABC}$ lying at $\ee^2$-distance from
the tangent plane $\pi_M$ to $\aa(\Omega_\ee)$ at $M$ 
is a curve lying in $\ee^2$-neighborhood of a quadric --- the
corresponding Dupin indicatrix of the surface $\aa(\Omega)$ in the 
central point $M$ of the $\ee$-patch.
\end{lemma}
\textsl{Proof.} 
Taking the appropriate Cartesian coordinate system
with the origin $M$ we can approximate the chosen $\ee$-patch by the
osculating paraboloid ${\cal P}=\{z=(K_1x^2+K_2y^2)/2\}$ with error terms of
order $\vrln{\vrln{o}}(\ee^2)$. 
We obtain the quadric $Q_M={\cal P} \cap \pi_{ABC}$ called the
Dupin indicatrix of $\aa(\Omega)$ at $M$. As one can easily check, the 
angle between $\pi_M$ and any of the tangent planes, taken at a
point of $\aa(\Omega_\ee)$, $\ee^2$-close to $Q_M$, is $\simeq \ee$:
the normals $\vec n =\aa_u \times \aa_v$ to $\aa(\Omega_\ee)$
at such points are $\ee^2$-close to the normals $\vec n_1$ of the
osculating paraboloid at the corresponding points with the same
$(x,y)$-coordinates. For the latter one has $\vec n_1(x,y)=-K_1x\vec
i + K_2 y \vec j + \vec k$ and for $x \simeq \ee$ and/or $y \simeq
\ee$, the angle between $\vec n_1(x,y)$ and $\vec n_1(0,0)=\vec k$
is $\simeq \ee$.

Now using  the fact that $\aa(\Omega_\ee)$ is $\ee^3$-close to the
osculating paraboloid and standard estimates for the values of
implicit functions and their derivatives we conclude that $I_{ABC}$
and $Q_M$ are $\ee^2$-close and the tangent directions at their
$\ee^2$-close points are also $\ee^2$-close. $\Box$

\end{document}